\documentclass[]{article}

\usepackage{authblk}
\usepackage{mathptmx}       % selects Times Roman as basic font
\usepackage{helvet}         % selects Helvetica as sans-serif font
\usepackage{courier}        % selects Courier as typewriter font
\usepackage{type1cm}        % activate if the above 3 fonts are
\usepackage{makeidx}         % allows index generation
\usepackage{graphicx}        % standard LaTeX graphics tool
\usepackage{multicol}        % used for the two-column index
\usepackage[bottom]{footmisc}% places footnotes at page bottom

\usepackage{psfrag}
\usepackage{amsfonts,amsmath,epsfig,amssymb,color,mathrsfs}
\newtheorem{teor}{Theorem}[section]
\newtheorem{lemma}[teor]{Lemma}
\newtheorem{defi}[teor]{Definition}
\newtheorem{remark}[teor]{Remark}
\newtheorem{ipot}[teor]{Assumption}
\newtheorem{prop}[teor]{Proposition}

\makeindex             % used for the subject index
\title{Conditional stability for backward parabolic operators with Osgood continuous coefficients}

\author[1]{Daniele Casagrande}
\author[2]{Daniele Del Santo}
\author[2]{Martino Prizzi}
\affil[1]{Dipartimento Politecnico di Ingegneria e Architettura, Universit\`a degli Studi di Udine, Via delle Scienze, 206 - 33100 Udine, Italy, \texttt{daniele.casagrande@uniud.it}}
\affil[2]{Dipartimento di Matematica e Geoscienze, Universit\`a degli Studi di Trieste, Via A. Valerio, 10 - 34100 Trieste, Italy, \texttt{delsanto@units.it,mprizzi@units.it}}

\begin{document}

\maketitle

\begin{abstract}
We prove continuous dependence on initial data for a backward parabolic operator whose leading coefficients are Osgodd continuous in time. This result fills the gap between uniqueness and continuity results obtained so far.
%
%
%
%The interest of the scientific community for the existence, uniqueness and stability of solutions to PDE's is testified by the numerous works available in the literature. In particular, in some recent publications on the subject~\cite{NONLINAL,MATAN} an inequality guaranteeing stability is shown to hold provided that the coefficients of the principal part of the differential operator are Log-Lipschitz continuous. Herein this result is improved along two directions. First, we describe how to construct an operator, whose coefficients in the principal part are not Log-Lipschitz continuous, for which the above mentioned inequality does not hold. Second, we show that the stability of the solution is guaranteed, in a suitable functional space, if the coefficients of the principal part are Osgood continuous.
\end{abstract}

\section{Introduction}
Backward parabolic equations are known to generate ill-posed (in the sense of Hadamard~\cite{Had_53,Had_64}) Cauchy problems. Due to the smoothing effects of the parabolic operator, in fact, it is not possible, in general, to guarantee existence of the solution for initial data which are not suitably regular. In addition, even when solutions possibly exist, uniqueness cannot be ensured without additional assumptions on the operator. Nevertheless, also for problems which are not well-posed the study of the conditional stability of the solution -- the surrogate of the notion of ``continuous dependence'' when existence of a solution is not guaranteed -- is of some interest. Such kind of study can be performed by resorting to the notion of \emph{well-behaved problem} introduced by John~\cite{John}: a problem is \emph{well-behaved} if ``only a fixed percentage of the significant digits need be lost in determining the solution from the data''. More precisely, a problem is well behaved if its solutions in a space $\mathcal{H}$ depend continuously on the data belonging to a space $\mathcal{K}$, provided they satisfy a prescribed bound in a space $\mathcal{H}^\prime$ (possibly different from $\mathcal{H}$).
In this paper we give a contribution to the study of the (\emph{well}) behaviour of the Cauchy problem associated with a backward parabolic operator. In particular, we consider the operator $\mathcal{L}$ defined, on the strip $[0,T]\times \mathbb{R}^n$, by 
\begin{equation}\label{eq_L}
\mathcal{L}u=\partial_t u +\sum_{i,j=1}^n \partial_{x_i}\left( a_{i,j}(t)\partial_{x_j}u \right)+\sum_{j=1}^n b_j(t)\partial_{x_j}u+c(t)u\,,
\end{equation}
where all the coefficients are bounded. We suppose that $a_{i,j}(t)=a_{j,i}(t)$ for all $i,j=1,\ldots,n$ and for all $t\in [0,T]$. We also suppose that $\mathcal{L}$ is backward parabolic, i.e. there exists $k_A\in ]0,1[$ such that, for all $(t,\xi)\in[0,T]\times \mathbb{R}^n$,
\begin{equation}
k_A\vert \xi\vert^2\le \sum_{i,j=1}^n a_{i,j}(t)\xi_i\xi_j\le k_A^{-1}\vert\xi\vert^2\,.
\end{equation}

We show that if the coefficients of the principal part of $\mathcal{L}$ are at least Osgood regular, then there exists a function space in which the associated Cauchy problem 
\begin{equation}\label{eq_probl_cauchy}
\left\{\begin{array}{ll}
\mathcal{L}u=f\,,\qquad & \textnormal{ in }(0,T)\times \mathbb{R}^n\,,\\
u\vert_{t=0}=u_0\,,\qquad & \textnormal{ in }\mathbb{R}^n\,,
\end{array}\right.
\end{equation}
has a stability property.

To collocate the new result in the framework of the existing literature, the contents of some publications on the subject are preliminarily recalled. They show that, as one could expect, the function space in which the stability property holds is related to the degree of regularity of the coefficients of $\mathcal{L}$. Weaker requirements on the regularity of the coefficients must be balanced, for the stability property to hold, by stronger \emph{a priori} requirements on the regularity of the solution, hence stability holds in a smaller function space.

The overview on available works helps to lead the reader to the new result, claimed in the final part of the paper, concerning operators with Osgood-continuous coefficients. This kind of regularity is critical since it is the minimum required regularity to have uniqueness of the solution and can therefore be considered as a sort of lower limit. Although the proof of the claim is based on the theoretical scheme followed to achieve previous results~\cite{MATAN}, the modifications needed to obtain an analogous proof in the case of Osgood coefficients are by no means trivial.

The paper is organised as follows. In Section~\ref{sec_uniq} we give an overview on uniqueness and non-uniqueness results for (\ref{eq_probl_cauchy}). Moreover, we introduce the notion of modulus of continuity and define the Osgood condition. Section~\ref{sec_stab} is dedicated to the notion of conditional stability; after recalling some known results, we state and prove the main result of the paper (Theorem~\ref{teo_nuovo}). In Section~\ref{sec_part} we consider the particular case of Log-Log-Lipschitz coefficients, where the dependence on initial data can be explicitly determined.

\section{Uniqueness and non-uniqueness results}\label{sec_uniq}
This section recalls some results on the uniqueness and non-uniqueness of the solution of the problem (\ref{eq_probl_cauchy}) for an operator like (\ref{eq_L}) with coefficients depending also on $x$. Consider the space
\begin{equation}\label{eq_insieme}
\mathcal{H}_0\triangleq C([0,T],L^2)\cap C([0,T[,H^1)\cap C^1([0,T[,L^2)\,.
\end{equation}
One of the first results concerning uniqueness is due to Lions and Malgrange~\cite{Lio_Mal} who consider an equation associated to a sesquilinear operator defined in a Hilbert space. In our context, this result can be read as follows.

\begin{teor}\label{teo_Lio_Mal}
If the coefficients of the principal part of $\mathcal{L}$ are Lipschitz continuous with respect to $t$ and $x$, $u\in\mathcal{H}_0$ and $u_0=0$, then $\mathcal{L}u=0$ implies $u\equiv 0$.$\hfill\square$
\end{teor}
The Lipschitz continuity of the coefficients is a crucial requirement for the claim, as shown some years later by Pli\'s~\cite{Pli} who proved the following theorem.

\begin{teor}
There exist $u$, $b_1$, $b_2$ and $c\in C^\infty(\mathbb{R}^3)$, bounded with bounded derivatives and periodic in the space variables and there exist $l:[0,T]\to\mathbb{R}$, H\"older-continuous of order $\delta$ for all $\delta<1$ but not Lipschitz-continuous, such that $1/2\le l(t)\le 3/2$ for all $t$, the support of $u$ is the set $\{t\ge 0\}\times\mathbb{R}^2$, and
\begin{multline}\label{eq_cauchy_plis}
%\left\{\begin{array}{ll}
\partial_t^2u(t,x_1,x_2)+\partial_{x_1}^2u(t,x_1,x_2)+l(t)\partial_{x_2}^2u(t,x_1,x_2)+\\
\qquad +b_1(t,x_1,x_2)\partial_{x_1}u(t,x_1,x_2)+b_2(t,x_1,x_2)\partial_{x_2}u(t,x_1,x_2)+\\
\qquad\qquad +c(t,x_1,x_2)u(t,x_1,x_2)=0 \qquad\qquad \textnormal{in }\mathbb{R}^3\,.%\\[2mm]
%\left.u(t,x_1,x_2)\right\vert_{t=0}=0 & \textnormal{in }\mathbb{R}
%\end{array}\right.\!\!\!
\end{multline}
%is the set $\mathbb{R}\times\mathbb{R}\times\{t\ge 0\}$.
$\hfill\square$
\end{teor}
Note that the differential operator in (\ref{eq_cauchy_plis}) is elliptic. However, the same idea developed by Pli\'s to prove the claim can be exploited to obtain a counterexample for the backward parabolic operator
$$
\mathcal{L}_P\triangleq\partial_t+\partial_{x_1}^2+l(t)\partial_{x_2}^2+b_1(t,x_1,x_2)\partial_{x_1}+b_2(t,x_1,x_2)\partial_{x_2}+c(t,x_1,x_2)\,.
$$
Moreover, the result can be extended to the operator $\mathcal{L}$ by considering the problem solved by $u(t,x_1,x_2)e^{-x_1^2-x_2^2}$, thus obtaining the following theorem.

\begin{teor}
There exist coefficients $a_{i,j}$, depending only on $t$, which are H\"older continuous of every order but not Lipschitz continuous and there exist $u\in\mathcal{H}_0$ such that the solution of problem (\ref{eq_probl_cauchy}) with $u_0=0$ and $f=0$ is not identically zero.$\hfill\square$
\end{teor}

In view of the previous results, a question naturally arises: which is the \emph{minimal} regularity (between Lipschitz continuity and H\"older continuity) of the coefficients of the principal part of $\mathcal{L}$ guaranteeing uniqueness of the solution of (\ref{eq_probl_cauchy})? To answer to this question, the definition of \emph{modulus of continuity}, that can be exploited to measure the degree of regularity of a function, is useful.

\begin{defi}
A \emph{modulus of continuity} is a function $\mu:[0,1]\to[0,1]$ which is continuous, increasing, concave and such that $\mu(0)=0$. A function $f:\mathbb{R}\to\mathbb{R}$ has \emph{regularity $\mu$} if
$$
\sup_{0<\vert t-s\vert<1}\frac{\vert f(t)-f(s)\vert}{\mu(\vert t-s\vert)}<+\infty\,.
$$
The set of all functions having regularity $\mu$ is denoted by $C^\mu$.
\end{defi}

As particular cases, the Lipschitz continuity, the $\tau$-H\"older continuity ($\tau\in]0,1[$) and the \emph{logarithmic Lipschitz} (in short \emph{Log-Lipschitz}) continuity are obtained for $\mu(s)=s$, $\mu(s)=s^{\tau}$ and $\mu(s)=s\log(1+1/s)$, respectively.

A further characterization of the modulus of continuity is the so called \emph{Osgood condition} which is crucial in most of the results on uniqueness and stability that are described in the rest of the article. A modulus of continuity $\mu$ satisfies the Osgood condition if
$$
\int_0^1 \frac{1}{\mu(s)}ds = +\infty\,.
$$

This characterization is used, for instance, in~\cite{JMPA} to obtain the following result concerning an operator whose coefficients in the principal part depend also on $x$.

\begin{teor}\label{teo_delpri_2005}
Let $\mu$ be a modulus of continuity that satisfies the Osgood condition. Let
\begin{equation}\label{eq_H_1}
\mathcal{H}_1\triangleq H^1([0, T],L^2(\mathbb{R}^n))\cap L^2([0, T],H^2(\mathbb{R}^n))
\end{equation}
and let the coefficients $a_{i,j}$ be such that, for all $i,j=1,\ldots,n$,
$$
a_{i,j}\in C^\mu([0,T],\mathscr{C}_b(\mathbb{R}^n))\cap \mathscr{C}([0,T],\mathscr{C}_b^2(\mathbb{R}^n))\,,
$$
where $\mathscr{C}_b$ is the space of bounded functions and $\mathscr{C}_b^2$ is the space of the bounded functions whose first and second derivatives are bounded. If $u\in \mathcal{H}_1$, $\mathcal{L}u=0$ on $[0,T]\times \mathbb{R}^n$ and $u(0,x)=0$ on $\mathbb{R}^n$, then $u\equiv 0$ on $[0,T]\times\mathbb{R}^n$.
\end{teor}

More recently, by using Bony's para-multiplication, the result has been improved as far as the regularity with respect to $x$ is concerned, i.e. replacing $\mathscr{C}^2$ regularity with Lipschitz regularity~\cite{AMPA}.

Note that the claim of Theorem~\ref{teo_delpri_2005} refers to the function space defined by (\ref{eq_H_1}), however, it is not difficult to extend it to the function space $\mathcal{H}_0$ defined by (\ref{eq_insieme}).

\section{Conditional stability results}\label{sec_stab}
For Cauchy problems related to the backward parabolic differential operators, which in general are not well posed, the notion of continuous dependence from initial data is replaced by the notion of (conditional) stability which is associated with the property of a problem to be well behaved, as defined by John~\cite{John}. The question about the conditional stability can be stated as follows. Suppose that two functions $u$ and $v$, defined in $[0,T]\times\mathbb{R}^n$, are solutions of the same equation; suppose, in addition, that $u$ and $v$ satisfy a fixed bound in a space $\mathcal{K}$ and that $\Vert u(0,\cdot)-v(0,\cdot)\Vert_{\mathcal{H}}$ is small (less than some $\epsilon$). Given these assumptions can we say something on the quantity $\sup_{t\in[0,T^\prime]}\Vert u(t,\cdot)-v(t,\cdot)\Vert_{\mathcal{K}}$ for some $T^\prime<T$? Does it remains small as well (e.g. less than a value related to $\epsilon$)? In this section some results that give an answer to the above questions are reported.

\subsection{Stability with Lipschitz-continuous (with respect to $t$) coefficients}
One of the first results on conditional stability has been proven by Hurd~\cite{Hurd} in the same theoretical framework considered by Lions and Malgrange.

\begin{teor}\label{teor_LioMal}
Suppose that the coefficients $a_{i,j}$ are Lipschitz continuous both in $t$ and in $x$. For every $T^\prime\in]0,T[$ and for every $D>0$ there exist $\rho>0$, $\delta\in]0,1[$ and $M>0$ such that if $u\in \mathcal{H}_0$ is a solution of $\mathcal{L}u=0$ on $[0,T]\times\mathbb{R}^n$ with $\Vert u(t,\cdot)\Vert_{L^2}\le D$ on $[0,T]$ and $\Vert u(0,\cdot)\Vert_{L^2}\le\rho$, then
\begin{equation}\label{eq_teo_hurd}
\sup_{t\in[0,T^\prime]}\Vert u(t,\cdot)\Vert_{L^2}\le M\Vert u(0,\cdot)\Vert_{L^2}^\delta\,.
\end{equation}
The constants $\rho$, $\delta$ and $M$ depend only on $T^\prime$ and $D$, on the ellipticity constant of $\mathcal{L}$, on the $L^\infty$ norms of the coefficients $a_{i,j}$, $b_j$, $c$, on the $L^\infty$ norms of their spatial derivatives, and on the Lipschitz constant of the coefficients $a_{i,j}$ with respect to time.$\hfill\square$
\end{teor}

The result expressed by (\ref{eq_teo_hurd}) implies uniqueness of the solution to the Cauchy problem, so that a necessary condition to this kind of conditional stability is that the coefficients $a_{i,j}$ fulfil the Osgood condition with respect to time. Hence a natural question arises: is Osgood condition also a sufficient condition for (\ref{eq_teo_hurd}) to hold? Del Santo and Prizzi~\cite{MATAN} have given a negative answer to this question. In particular, mimicking Pli\'s counterexample, they have shown that if the coefficients $a_{i,j}$ are not Lipschitz-continuous but only Log-Lipschitz-continuous then Hurd's result does not hold. Moreover, they have proven that if the coefficients are Log-Lipschitz-continuous then a conditional stability property, weaker than (\ref{eq_teo_hurd}), does hold. More recently, the result has been further improved~\cite{NONLINAL}.

\subsection{Stability with Log-Lipschitz-continuous (with respect to $t$) coefficients}
As mentioned above, Osgood condition is not sufficient for H\"older conditional stability of the solution expressed by (\ref{eq_teo_hurd}). The following paragraph specifies this claim.
\subsubsection{Counterexample to H\"older stability in the Log-Lipschitz case}\label{par_controesempio_uno}
The counterexample relies on the fact that it is possible~\cite{MATAN} to construct
\begin{itemize}
\item a sequence $\{\mathcal{L}_k\}_{k\in \mathbb{N}}$ of backward uniformly parabolic operators with uniformly Log-Lipschitz-continuous coefficients (not depending on the space variables) in the principal part and space-periodic uniformly bounded smooth coefficients in the lower order terms,
\item a sequence $\{u_k\}_{k\in\mathbb{N}}$ of space-periodic smooth uniformly bounded solutions of $\mathcal{L}_ku_k=0$ on $[0,1]\times\mathbb{R}^2$,
\item a sequence $\{t_k\}_{k\in\mathbb{N}}$ of real numbers, with $t_k\to 0$,
\end{itemize}
such that
$$
\lim_{k\to\infty}\Vert u_k(0,\cdot,\cdot)\Vert_{L^2([0,2\pi]\times[0,2\pi])}=0
$$
and
$$
\lim_{k\to\infty}\frac{\Vert u_k(t_k,\cdot,\cdot)\Vert_{L^2([0,2\pi]\times[0,2\pi])}}{\Vert u_k(0,\cdot,\cdot)\Vert^\delta_{L^2([0,2\pi]\times[0,2\pi])}}=+\infty
$$
for every $\delta>0$.
%We remark that this situation is exactly what is needed to show that for backward operators with Log-Lipscitz continuous coefficient a result similar to Theorem~\ref{teor_LioMal} cannot hold.

\subsubsection{Stability result in the Log-Lipschitz case}
In the case of Log-Lipschitz coefficients a result weaker that (\ref{eq_teo_hurd}) is valid. Consider the equation $\mathcal{L}u=0$ on $[0,T]\times\mathbb{R}^n$, with $\mathscr{L}$ defined in (\ref{eq_L}) and suppose that for all $i,j=1,\ldots,n$, $a_{i,j}\in \textnormal{LogLip}([0,T])$, in particular
$$
\sup_{0<\vert\tau\vert<1}\frac{\vert a_{i,j}(t+\tau)-a_{i,j}(t)\vert}{\vert\tau\vert\left(\log\left(1+\frac{1}{\vert\tau\vert}\right)\right)}<+\infty\,;\label{condizione_terza}
$$
let $b_j$ and $c$ belong to $L^\infty([0,T])$.

\begin{teor}\label{teo_matan}
{\bf\cite{MATAN}}
Suppose that the above hypotheses hold. For all $T^\prime\in]0,T[$ and for all $D>0$ there exist $\rho>0$, $M>0$, $N>0$ and $\delta\in]0,1[$ such that, if $u\in\mathcal{H}_0$ is a solution of $\mathcal{L}u=0$ on $[0,T]\times\mathbb{R}^n$ with $\Vert u(t,\cdot)\Vert_{L^2}\le D
$ on $[0,T]$ and $\Vert u(0,\cdot)\Vert_{L^2}\le \rho$, then
\begin{equation}\label{eq_teo_del_pri}
\sup_{t\in[0,T^\prime]}\Vert u(t,\cdot)\Vert_{L^2}\le M e^{-N\vert\log\Vert u(0,\cdot)\Vert_{L^2}\vert^\delta}\,,
\end{equation}
where the constants $\rho$, $\delta$, $M$ and $N$ depend only on $T^\prime$, on $D$, on the ellipticity constant of $\mathcal{L}$, on the $L^\infty$ norms of the coefficients $a_{i,j}$, on the $L^\infty$ norms of their spatial first derivatives, and on the Log-Lipschitz constant of the coefficients $a_{i,j}$ with respect to time.
\end{teor}

Using Bony's para-product the result can be extended to the case in which the coefficients depend also on the space variable and are Lipschitz continuous with respect to it~\cite{NONLINAL}.

%\begin{teor}\label{teo_nonlinal}
%{\bf\cite{NONLINAL}}
%The same result holds if hypotheses \ref{condizione_terza}, \ref{condizione_quarta}, \ref{condizione_ultima} are replaced by
%\begin{itemize}
%\item[3')] for all $i,j=1,\ldots,n$, $a_{i,j}\in \textnormal{LogLip}([0,T],L^\infty(\mathbb{R}^n))\cap L^\infty([0,T],\textnormal{Lip}(\mathbb{R}^n))$;
%\item[4')] for all $j=1,\ldots,n$, $b_j\in L^\infty([0,T],\textnormal{Lip}(\mathbb{R}^n))$;
%\item[5')] $c\in L^\infty([0,T],\textnormal{Lip}(\mathbb{R}^n))$.$\hfill\square$
%\end{itemize}
%\end{teor}

\subsection{Stability with Osgood-continuous (with respect to time) coefficients}
Let us finally come to the new result contained in this paper. As in the previous section we first present a counterexample to the stability condition (\ref{eq_teo_del_pri}) and then a new weaker stability result.

\subsubsection{Counterexample to stability estimate (\ref{eq_teo_del_pri}) in the LogLog-Lipschitz case}\label{par_controesempio_due}
Consider the modulus of continuity $\omega$ defined, near $0$, by
$$
\omega(s)=s\log\left(1+\frac{1}{s}\right)\log\left(\log\left(1+\frac{1}{s}\right)\right)
$$
and note that $\omega$ satisfies the Osgood condition but $\mathscr{C}^\omega$ functions are not Log-Lipschitz continuous. As in Paragraph~\ref{par_controesempio_uno}, it is possible~\cite{mia_tesi} to construct
\begin{itemize}
\item a sequence $\{\mathcal{L}_k\}_{k\in \mathbb{N}}$ of backward uniformly parabolic operators with uniformly $\mathscr{C}^\omega$-continuous coefficients in the principal part and space-periodic uniformly bounded smooth coefficients in the lower order terms,
\item a sequence $\{u_k\}_{k\in\mathbb{N}}$ of space-periodic smooth uniformly bounded solutions of $\mathcal{L}_ku_k=0$ on $[0,1]\times\mathbb{R}^2$,
\item a sequence $\{t_k\}_{k\in\mathbb{N}}$ of real numbers, with $t_k\to 0$,
\end{itemize}
such that
$$
\lim_{k\to\infty}\Vert u_k(0,\cdot,\cdot)\Vert_{L^2([0,2\pi]\times[0,2\pi])}=0
$$
but (\ref{eq_teo_del_pri}) does not hold for all $\delta$; more precisely
$$
\lim_{k\to\infty}\frac{\Vert u_k(t_k,\cdot,\cdot)\Vert_{L^2([0,2\pi]\times[0,2\pi])}}{e^{-N\vert \log\Vert u_k(0,\cdot,\cdot)\Vert_{L^2([0,2\pi]\times[0,2\pi])}\vert^\delta}}=+\infty
$$
for every $\delta>0$.

\subsection{Stability result in the Osgood-continuous case}
%Consider the operator (\ref{eq_L}) and suppose that its coefficients only depend on time, namely consider the equation
%\begin{equation}\label{eq_retro}
%\partial_t{u}+\sum_{i,j=1}^na_{i,j}(t)\partial_{x_i,x_j}u+\sum_{i=1}^nb_i(t)\partial_{x_i}u+c(t)u=0
%\end{equation}
%on the strip $[0,T]\times \mathbb{R}^n$.
%\begin{defi}
%Consider the operator (\ref{eq_L}). A function $u$, belonging to the functional space
%\begin{equation}\label{eq_definizione_spazio}
%\mathcal{E}\triangleq \mathscr{C}^0([0,T],L^2(\mathbb{R}^n))\cap \mathscr{C}^0([0,T),H^1(\mathbb{R}^n))\cap\mathscr{C}^1([0,T),L^2(\mathbb{R}^n))
%\end{equation}
%is a \emph{solution} of
%\begin{equation}\label{eq_lu=0}
%\mathscr{L}u=0
%\end{equation}
%if it verifies (\ref{eq_lu=0}) in $H^{-1}(\mathbb{R}^n)$.$\hfill\Diamond$
%\end{defi}
From now on, the following conditions are assumed to hold.

\medskip

\begin{ipot}\label{ipot_princ}
The operator $\mathcal{L}$ defined in (\ref{eq_L}) is such that
\begin{itemize}
\item for all $t\in[0,T]$ and for all $i,j=1,\ldots,n$,
$$
a_{i,j}(t)=a_{j,i}(t)\,;
$$
\item there exists $k_A>0$ such that, for all $(t,\xi)\in [0,T]\times \mathbb{R}^n$,
$$
k_A\vert\xi\vert^2\le \sum_{i,j=1}^na_{i,j}(t)\xi_i\xi_j\le k_A^{-1}\vert\xi\vert^2\,;
$$
\item there exists $k_B>0$ such that, for all $t\in [0,T]$ and for all $i=1,\ldots,n$, $\vert b_i(t)\vert \le k_B$;
\item there exists $k_C>0$ such that, for all $t\in [0,T]$, $\vert c(t)\vert\le k_C$;
\item for all $i,j=1,\ldots,n$, $a_{i,j}\in C^{\omega}([0,T])$, where $\omega$ is a modulus of continuity that satisfies the Osgood condition.
\end{itemize}
\end{ipot}

\medskip

We can now state our main result.

\medskip

\begin{teor}\label{teo_nuovo}
For all $T^\prime\in]0,T[$ and for all $D>0$ there exist $\rho^\prime>0$, and an increasing continuous function $G:[0,+\infty[\to[0,+\infty[$, with $G(0)=0$, such that, if $u\in\mathcal{H}_0$ is a solution of $\mathcal{L}u=0$ on $[0,T]$ with $\Vert u(t,\cdot)\Vert_{L^2}\le D$ on $[0,T]$ and $\Vert u(0,\cdot)\Vert_{L^2}\le\rho^\prime$, then
\begin{equation}\label{eq_in_L2}
\sup_{t\in[0,T^\prime]}\Vert u(t,\cdot)\Vert_{L^2}^2\le G(\Vert u(0,\cdot)\Vert_{L^2}^2)\,.
\end{equation}
The constant $\rho^\prime$ and the function $G$ depend on $k_A,k_B,k_C,\omega,n,T,T^\prime$ and $D$.$\hfill\square$
\end{teor}

\begin{defi}\label{def_Gev_Sob}
\cite{Hua_Rod} Given $a\ge 0$, $d\in\mathbb{R}$ and $\epsilon>1$, the \emph{Gevrey-Sobolev} function space $H^d_{a,\epsilon}$ is the space of the functions $u:\mathbb{R}^n\to\mathbb{R}$ such that
$$
\Vert u\Vert_{H^d_{a,\epsilon}}\triangleq\int_{\mathbb{R}^n}\left(1+\vert\xi\vert^2\right)^de^{2a\vert\xi\vert^{1/\epsilon}}\left\vert \hat{u}(\xi)\right\vert^2d\xi<+\infty\,,
$$
where $\hat{u}$ is the Fourier transform of $u$.
\end{defi}

\begin{defi}\label{def_OGS}
Let $a>0$, $d\in\mathbb{R}$ and $\omega$ a modulus of continuity satisfying the Osgood condition. We denote by $H_{a,\omega}^d$ the set of the functions $u:\mathbb{R}^n\to\mathbb{R}$ such that
$$
\Vert u\Vert^2_{H_{a,\omega}^d}\triangleq \int_{\mathbb{R}^n}\left(1+\vert\xi\vert^2\right)^de^{a\vert\xi\vert^2\omega\left(\frac{1}{\vert\xi\vert^2+1}\right)}\vert\hat{u}(\xi)\vert^2d\xi<+\infty\,.
$$
We call it \emph{Osgood-Sobolev} function space.
\end{defi}

\begin{remark}
From Definitions~\ref{def_Gev_Sob} and~\ref{def_OGS} it is easy to see that, for all moduli of continuity $\omega$, for all $\epsilon>1$, for all $a>0$ and for all $d\in\mathbb{R}$,
$$
H_{a,\omega}^d \subset H_{a,\epsilon}^d\,.
$$
\end{remark}

Theorem~\ref{teo_nuovo} is a consequence of the following local result.

\begin{teor}\label{teo_finale_stab_norme_classiche}
There exists $\alpha_1>0$ and, for any $T^{\prime\prime}:0<T^{\prime\prime}<T$, there exist constants $\rho>0$, $C>0$ and a function $g:[0,k_A]\to[0,+\infty[$, such that, if $u\in\mathcal{H}_0$ is a solution of
\begin{equation}\label{eq_lu=0}
\mathcal{L}u=0\,,
\end{equation}
with $\mathcal{L}$ fulfilling Assumption~\ref{ipot_princ} and $\Vert u(0,\cdot)\Vert^2_{H_{\nu,\epsilon}^0}<\rho$ for some $\nu>0$ and some $\epsilon>1$, then
\begin{equation}\label{eq_u_settima_ultima_norme_classiche}
\sup_{z\in[0,\bar{\sigma}]}\Vert u(z,\cdot)\Vert^2_{H^1}\le C e^{-\sigma g\left(\Vert u(0,\cdot)\Vert^2_{H_{\nu,\epsilon}^0}\right)}\left[1+\Vert u(\sigma,\cdot)\Vert^2_{H^1}\right]\,,
\end{equation}
where $\sigma=\min\{T^{\prime\prime},1/\alpha_1\}$ and $\bar{\sigma}=\sigma/8$. The constant $\alpha_1$ depends only on $k_A,k_B,k_C,\omega$ and $n$ while the constants $\rho$ and $C$ depend also on $T$ and $T^{\prime\prime}$. The function $g$ is a strictly decreasing function; it depends on $k_A,k_B,k_C,\omega,n,T,T^{\prime\prime},\epsilon$ and $\nu$ and satisfies $\lim_{y\to 0}g(y)=+\infty$.$\hfill\square$
\end{teor}
Theorem~\ref{teo_finale_stab_norme_classiche} will be proven with the help of partial results expressed in terms of estimates of some integral quantities. The following Lemma~\ref{lem_tutti gevrey} guarantees that all the integral quantities that will be introduced are finite, so that the obtained estimates make sense.

\begin{lemma}\label{lem_derivata_limitata}
Let $u:[0,T]\to\mathbb{R}$ a $C^1$ function. If $u^\prime(t)\ge M u(t)$, then $u(t)\le e^{M(t-T)}u(T)$.
\end{lemma}
\textbf{Proof. }If is sufficient to note that:
\begin{multline*}
u^\prime(t)\ge Mu(t)\,\,\Rightarrow\,\,u^\prime(t)e^{-M(t-T)}-Mu(t)e^{-M(t-T)}\ge 0\,\,\Rightarrow\\[2mm]
\Rightarrow\,\,\frac{d}{dt}\left(u(t)e^{-M(t-T)}\right)\ge 0\,\,\Rightarrow\,\, u(t)e^{-M(t-T)}\le u(T)\,\,\Rightarrow\\
\Rightarrow\,\,u(t)\le e^{M(t-T)}u(T)\,.
\end{multline*}
$\hfill\blacksquare$

%Let $\hat{u}$ denote the Fourier transform of $u$ with respect to $x$ and let $H_{a,\epsilon}^d$ denote a Gevrey-Sobolev function space, i.e. the space of the functions $u:\mathbb{R}^n\to\mathbb{R}$ such that
%$$
%\int_{\mathbb{R}^n}\left(1+\vert\xi\vert^2\right)^de^{2a\vert\xi\vert^{1/epsilon}}\left\vert \hat{u}(\xi)\right\vert^2d\xi<+\infty\,.
%$$

\begin{lemma}\label{lem_tutti gevrey}
Let $M>0$ and let $u\in\mathcal{H}_0$ be a solution of
\begin{equation}\label{eq_lemma_maggiorazioni}
\partial_t u+\sum_{i,j=1}^na_{i,j}(t)\partial_{x_i}\partial_{x_j}u+\sum_{i=1}^nb_i(t)\partial_{x_i}u+c(t)u=0\,,
\end{equation}
on $[0,T]$, such that $\Vert u(t,\cdot)\Vert_{L^2}\le M$, for all $t\in[0,T]$. Let $l>0$ and extend the coefficients $a_{i,j}$, $b_i$ and $c$ to $[-l,T]$ by setting $a_{i,j}(t)=a_{i,j}(0)$, $b_i(t)=b_i(0)$ and $c(t)=c(0)$ for all $t\in[-l,0[$. Then $u$ can be extended to a solution of (\ref{eq_lemma_maggiorazioni}) on $[-l,T]$ such that there exists $\hat{M}$ such that $\Vert u(t,\cdot)\Vert_{L^2}\le \hat{M}$ on $[-l,T]$. The constant $\hat{M}$ depends only on $n$, $k_A$, $k_B$, $K_C$, $T$, $l$ and $M$. Moreover,
\begin{enumerate}
\item $u\in C^0([-l,T[,H_{a,\epsilon}^d)$ for all $a\ge 0$, $\epsilon>1$ and $d\in\mathbb{R}$;
\item $u\in C^0([-l,T[,H^1)$ and there exists $C$, which depends on $n$, $k_A$, $k_B$, $k_C$, $T$ and $l$, such that
$$
\Vert u(t,\cdot)\Vert_{H^1}\le C(T-t)^{-1/2}\Vert u(T,\cdot)\Vert_{L^2}
$$
for all $t\in[-l,T[$;
\item there exists $\hat{C}$, which depends on $n$, $k_A$, $k_B$, $k_C$, $l$, $\nu$ and $\epsilon$ and which tends to $+\infty$ when $l$ tends to zero, such that
$$
\Vert u(-l,\cdot)\Vert_{H^0_{\nu,\epsilon}}\le \hat{C}\Vert u(0,\cdot)\Vert_{L^2}\,.
$$
$\hfill\square$
\end{enumerate}
\end{lemma}
\textbf{Proof. }It is easy to see that for all $t\in[0,T]$ and for almost all $\xi\in\mathbb{R}^n$,
\begin{equation}\label{identita_trasf_four}
\partial_t \hat{u}(t,\xi)-\sum_{i,j=1}^n a_{i,j}(t)\xi_i\xi_j \hat{u}(t,\xi)+\imath\sum_{i=1}^nb_i(t)\xi_i\hat{u}(t,\xi)+c(t)\hat{u}(t,\xi)=0\,.
\end{equation}
Multiplying both terms of (\ref{identita_trasf_four}) by $\bar{\hat{u}}$ yields
\begin{equation}\label{eq_ut_e_u_bar}
\partial_t \hat{u}(t,\xi)\bar{\hat{u}}(t,\xi)\!=\!\sum_{i,j=1}^n a_{i,j}(t)\xi_i\xi_j \vert\hat{u}(t,\xi)\vert^2\!-\!\imath\sum_{i=1}^nb_i(t)\xi_i\vert\hat{u}(t,\xi)\vert^2\!-\!c(t)\vert\hat{u}(t,\xi)\vert^2\,.
\end{equation}
By adding to (\ref{eq_ut_e_u_bar}) its complex conjugate, we obtain
\begin{multline}\label{eq_u_tutti_una}
\partial_t\vert\hat{u}(t,\xi)\vert^2=2\sum_{i,j=1}^na_{i,j}(t)\xi_i\xi_j\vert\hat{u}(t,\xi)\vert^2+2\sum_{i=1}^n\Im\{b_i(t)\}\xi_i\vert \hat{u}(t,\xi)\vert^2+\\[2mm]
-2\Re\{c(t)\}\vert\hat{u}(t,\xi)\vert^2\,,
\end{multline}
hence, recalling the bounds for the coefficients of $\mathcal{L}$ (see Assumption~\ref{ipot_princ}), 
\begin{equation*}%\label{eq_u_tutti_una}
\partial_t\vert\hat{u}(t,\xi)\vert^2\ge2 k_A\vert\xi\vert^2\vert\hat{u}(t,\xi)\vert^2-2n k_B\vert\xi\vert\vert\hat{u}(t,\xi)\vert^2-2k_C\vert\hat{u}(t,\xi)\vert^2\,,
\end{equation*}
i.e.
\begin{equation*}%\label{eq_u_tutti_due}
\partial_t\vert\hat{u}(t,\xi)\vert^2\ge(2k_A\vert\xi\vert^2-2n k_B\vert\xi\vert-2k_C)\vert\hat{u}(t,\xi)\vert^2\,.
\end{equation*}
Lemma~\ref{lem_derivata_limitata} allows one to write
\begin{equation}\label{tre_asterischi}
\vert\hat{u}(t,\xi)\vert^2\le e^{(2k_A\vert\xi\vert^2-2nk_B\vert\xi\vert-2k_C)(t-T)}\vert \hat{u}(T,\xi)\vert^2\,.
\end{equation}
Therefore, for a fixed $t\in[-l,T[$,
\begin{multline*}
\int_{\mathbb{R}^n}\left(1+\vert\xi\vert^2\right)^{d}e^{2a\vert\xi\vert^{\frac{1}{\epsilon}}}\vert\hat{u}(t,\xi)\vert^2d\xi\le\\
\le\int_{\mathbb{R}^n}\left(1+\vert\xi\vert^2\right)^{d}e^{2a\vert\xi\vert^{\frac{1}{\epsilon}}+(2k_A\vert\xi\vert^2-2nk_B\vert\xi\vert-2k_C)(t-T)}\vert\hat{u}(T,\xi)\vert^2d\xi<+\infty\,,
\end{multline*}
where the last inequality comes from the fact that $u\in\mathcal{H}_0$ and therefore, in particular, $u\in\mathscr{C}^0([0,T],L^2(\mathbb{R}^n))$, and, since $t<T$,
$$
\lim_{\vert \xi \vert\to\infty}\left(1+\vert\xi\vert^2\right)^{d}e^{2a\vert\xi\vert^{\frac{1}{\epsilon}}+(2k_A\vert\xi\vert^2-2nk_B\vert\xi\vert-2k_C)(t-T)}=0
$$
for all $a>0$ and all $\epsilon>1$. The first claim is then proven. The second claim is proven easily by choosing $d=1$ and $a=0$. To prove the third claim it is sufficient to rewrite equation (\ref{tre_asterischi}) replacing $T$ with $0$.$\hfill\blacksquare$

\subsection{Preliminary results and defintions}\label{sec_omega}
In this section some functions that are used in the rest of the article are defined. Let $\omega$ be a modulus of continuity satisfying Osgood condition. For a given $\rho>1$ define the function $\theta:[1,+\infty[\to[0,+\infty]$ as
\begin{equation}\label{eq_def_theta}
\theta(\rho)=\int_{1/\rho}^1\frac{1}{\omega(s)}ds\,.
\end{equation}
It is easy to see that $\theta$ is bijective and strictly increasing. As a consequence, it can be inverted. For $y\in(0,1]$, for $q>0$ and for $\lambda>0$, let $\psi_{\lambda,q}:]0,1]\to[1,+\infty[$ be defined by

$$
\psi_{\lambda,q}(y)\triangleq \theta^{-1}\left(-\lambda q \log y\right)\,.
$$
The relation
$$
\theta\left(\psi_{\lambda,q}(y)\right)=-\lambda q \log y
$$
immediately follows from the definitions; hence
$$
\theta^\prime\left(\psi_{\lambda,q}(y)\right)\psi^\prime_{\lambda,q}(y)=-\frac{\lambda q}{y}\,.
$$
Now, let the function $\phi_{\lambda, q}:(0,1]\to(-\infty,0]$ be defined as
\begin{equation}\label{eq_def_phi}
\phi_{\lambda,q}(y)\triangleq q\int_1^y\psi_{\lambda,q}(z)dz\,.
\end{equation}
The function $\phi_{\lambda, q}$ is bijective and strictly increasing; moreover,
\begin{equation}\label{eq_diff_phi_prima}
\phi_{\lambda,q}^{\prime\prime}(y)=q\psi^\prime_{\lambda,q}(y)=\frac{q}{\theta^\prime\left(\psi_{\lambda,q}(y)\right)}\left(-\frac{\lambda q}{y}\right)\,.
\end{equation}
On the other hand, equation (\ref{eq_def_theta}), with the change of variable $\eta=1/s$, becomes
$$
\theta(\rho)=-\int_\rho^1\frac{1}{\omega\left(\frac{1}{\eta}\right)}\frac{1}{\eta^2}d\eta=\int_1^\rho\frac{1}{\eta^2\omega\left(\frac{1}{\eta}\right)}d\eta
$$
from which
\begin{equation}\label{eq_da_sostituire}
\frac{1}{\theta^\prime\left(\psi_{\lambda,q}(y)\right)}=\psi_{\lambda,q}(y)^2\omega\left(\frac{1}{\psi_{\lambda,q}(y)}\right)\,.
\end{equation}
Substituting (\ref{eq_da_sostituire}) into (\ref{eq_diff_phi_prima}) and recalling that $\psi_{\lambda,q}(y)=\phi^\prime_{\lambda,q}(y)/q$, it is easy to see that $\phi_{\lambda,q}$ satisfies the equation
\begin{equation}\label{eq_diff_per_phi}
y\phi^{\prime\prime}_{\lambda,q}(y)=-\lambda\left(\phi^\prime_{\lambda,q}(y)\right)^2\omega\left(\frac{q}{\phi_{\lambda,q} ^\prime(y)}\right)\,.
\end{equation}
Note that for all $\lambda>0$, for all $q>0$ and for all $y\in(0,1]$, $\psi_{\lambda,q}\in(1,+\infty)$ and, consequently,
$$
\frac{q}{\phi^\prime_{\lambda,q}(y)}\in(0,1)\,.
$$

\medskip

\subsection{A pointwise estimate}
The first result shows that, once fixed $\xi$, namely the value of the frequence argument of $\hat{u}$, it is possible to find a bound for a particular time-integral, in an interval $[0,\sigma]$, of a function of $\vert\hat{u}(t,\xi)\vert$. This bound consists in the sum of two terms depending on $\hat{u}(0,\xi)$ and $\hat{u}(\sigma,\xi)$, respectively.

\begin{prop}\label{prop_con_lettere}
Let $T^{\prime\prime}\in \; ]0,T[$. There exist $\alpha_1>0$, $\bar{\lambda}$ and $\bar{\gamma}>0$ such that, setting $\alpha\triangleq\max\{\alpha_1,1/T^{\prime\prime}\}$, defining $\sigma=1/\alpha$, fixing $\tau\in\;]0,\sigma/4]$, and letting $\beta\ge \sigma+\tau$, whenever $u\in \mathcal{H}_0$ is a solution of (\ref{eq_lu=0}), one has
\begin{multline}\label{eq_u_hat_seconda}
\frac{1}{4}\left(k_A\vert\xi\vert^2+\gamma\right)\int_0^\sigma e^{(1-\alpha t)\vert\xi\vert^2\omega\left(\frac{1}{\vert\xi\vert^2+1}\right)}e^{2\gamma t}e^{-2\beta\phi_\lambda\left(\frac{t+\tau}{\beta}\right)}\vert\hat{u}(t,\xi)\vert^2dt\le\\
\le \phi^\prime_\lambda\left(\frac{\tau}{\beta}\right)\tau e^{\vert\xi\vert^2\omega\left(\frac{1}{\vert\xi\vert^2+1}\right)}e^{-2\beta\phi_\lambda\left(\frac{\tau}{\beta}\right)}\vert\hat{u}(0,\xi)\vert^2+\\
+(\sigma+\tau)(\gamma+k_A^{-1}\vert\xi\vert^2)e^{2\gamma\sigma}e^{-2\beta\phi_\lambda\left(\frac{\sigma+\tau}{\beta}\right)}\vert\hat{u}(\sigma,\xi)\vert^2\,,
\end{multline}
for all $\lambda\ge\bar{\lambda}$ and all $\gamma\ge\bar{\gamma}$, where $\phi_\lambda\triangleq \phi_{\lambda,k_A}$ (see (\ref{eq_def_phi})).
%\begin{equation}\label{eq_definzione_phi}
%y\,\phi^{\prime\prime}(y)=-\lambda(\phi^\prime(y))^2\omega\left(\frac{k_A}{\phi^\prime(y)}\right)
%\end{equation}
The constant $\alpha_1$ depends only on $n$, $k_A$, $k_B$, $k_C$ and $\omega$, while $\bar{\gamma}$ and $\bar{\lambda}$ depend on $n$, $k_A$, $k_B$, $k_C$, $\omega$, $T$ and $T^{\prime\prime}$.$\hfill\square$
\end{prop}

%There exist $\alpha_1>0$, $\bar{\lambda}$ and $\bar{\gamma}>0$ such that, setting $\alpha\triangleq\max\{\alpha_1,1/T^{\prime\prime}\}$ for some $T^{\prime\prime}:0<T^{\prime\prime}<T$, defining $\sigma=1/\alpha$, fixing $\tau\in]0,\sigma/4]$, and letting $\beta\ge \sigma+\tau$ be a free parameter, whenever $u\in \mathcal{H}_0$ is a solution of (\ref{eq_lu=0}), one has
%\begin{multline}\label{eq_u_hat_seconda}
%\frac{1}{4}\left(k_A\vert\xi\vert^2+\gamma\right)\int_0^\sigma e^{(1-\alpha t)\vert\xi\vert^2\omega\left(\frac{1}{\vert\xi\vert^2+1}\right)}e^{2\gamma t}e^{-2\beta\phi_\lambda\left(\frac{t+\tau}{\beta}\right)}\vert\hat{u}(t,\xi)\vert^2dt\le\\\le \phi^\prime_\lambda\left(\frac{\tau}{\beta}\right)\tau e^{\vert\xi\vert^2\omega\left(\frac{1}{\vert\xi\vert^2+1}\right)}e^{-2\beta\phi_\lambda\left(\frac{\tau}{\beta}\right)}\vert\hat{u}(0,\xi)\vert^2+\\+(\sigma+\tau)(\gamma+k_A^{-1}\vert\xi\vert^2)e^{2\gamma\sigma}e^{-2\beta\phi_\lambda\left(\frac{\sigma+\tau}{\beta}\right)}\vert\hat{u}(\sigma,\xi)\vert^2\,,\end{multline}for all $\lambda\ge\bar{\lambda}$ and all $\gamma\ge\bar{\gamma}$, where $\phi_\lambda\triangleq \phi_{\lambda,k_A}$ (see (\ref{eq_diff_per_phi})).
%\begin{equation}\label{eq_definzione_phi}
%y\,\phi^{\prime\prime}(y)=-\lambda(\phi^\prime(y))^2\omega\left(\frac{k_A}{\phi^\prime(y)}\right)
%\end{equation} The constant $\alpha_1$ depends only on $n$, $k_A$, $k_B$, $k_C$ and $\omega$, while $\bar{\gamma}$ and $\bar{\lambda}$ depend on $n$, $k_A$, $k_B$, $k_C$, $\omega$, $T$ and $T^{\prime\prime}$.$\hfill\square$
%\end{prop}

\textbf{Proof. } Let $T^{\prime\prime}\in \; ]0,T[$ and let $\alpha\geq1/T^{\prime\prime}$, $\gamma>0$, $\lambda>0$, $\tau\in\;]0,T^{\prime\prime}[$, $\sigma=1/\alpha$ and $\beta\geq \tau+\sigma$. Consider the function $\hat{v}:[0,\sigma]\times\mathbb{R}^n\to\mathbb{R}$ defined by
\begin{equation}\label{eq_v_hat}
\hat{v}(t,\xi)=e^{\left(\frac{1-\alpha t}{2}\right)\vert\xi\vert^2\omega\left(\frac{1}{\vert\xi\vert^2+1}\right)}e^{\gamma t}e^{-\beta\phi_\lambda\left(\frac{t+\tau}{\beta}\right)}\hat{u}(t,\xi)\,.
\end{equation}
The time-derivative of $\hat{v}$ is
\begin{multline*}
\partial_t\hat{v}(t,\xi)=-\frac{\alpha}{2}\vert\xi\vert^2\omega\left(\frac{1}{\vert\xi\vert^2+1}\right)e^{\left(\frac{1-\alpha t}{2}\right)\vert\xi\vert^2\omega\left(\frac{1}{\vert\xi^2\vert+1}\right)}e^{\gamma t}e^{-\beta\phi_\lambda\left(\frac{t+\tau}{\beta}\right)}\hat{u}(t,\xi)+\\
+\gamma e^{\left(\frac{1-\alpha t}{2}\right)\vert\xi\vert^2\omega\left(\frac{1}{\vert\xi\vert^2+1}\right)}e^{\gamma t}e^{-\beta\phi_\lambda\left(\frac{t+\tau}{\beta}\right)}\hat{u}(t,\xi)+\\
-\phi_\lambda^\prime\left(\frac{t+\tau}{\beta}\right)e^{\left(\frac{1-\alpha t}{2}\right)\vert\xi\vert^2\omega\left(\frac{1}{\vert\xi\vert^2+1}\right)}e^{\gamma t}e^{-\beta\phi_\lambda\left(\frac{t+\tau}{\beta}\right)}\hat{u}(t,\xi)+\\
+e^{\left(\frac{1-\alpha t}{2}\right)\vert\xi\vert^2\omega\left(\frac{1}{\vert\xi\vert^2+1}\right)}e^{\gamma t}e^{-\beta\phi_\lambda\left(\frac{t+\tau}{\beta}\right)}\partial_t\hat{u}(t,\xi)
\end{multline*}
which may be rewritten as
\begin{multline}\label{eq_per_v_hat}
\partial_t\hat{v}+\frac{\alpha}{2}\vert\xi\vert^2\omega\left(\frac{1}{\vert\xi\vert^2+1}\right)\hat{v}-\gamma\hat{v}+\phi_\lambda^\prime\left(\frac{t+\tau}{\beta}\right)\hat{v}-\sum_{i,j=1}^na_{i,j}(t)\xi_i\xi_j\hat{v}+\\
+\imath\sum_{i=1}^nb_i(t)\xi_i\hat{v}+c(t)\hat{v}=0\,,
\end{multline}
where the dependency of $\hat{v}$ and $\partial_t\hat{v}$ on $t$ and on $\xi$ has been neglected for the sake of a simple notation and where the identity (\ref{identita_trasf_four}) has been exploited. The complex conjugate equation of (\ref{eq_per_v_hat}) is
\begin{multline}\label{eq_per_v_hat_con}
\partial_t\bar{\hat{v}}+\frac{\alpha}{2}\vert\xi\vert^2\omega\left(\frac{1}{\vert\xi\vert^2+1}\right)\bar{\hat{v}}-\gamma\bar{\hat{v}}+\phi_\lambda^\prime\left(\frac{t+\tau}{\beta}\right)\bar{\hat{v}}-\sum_{i,j=1}^na_{i,j}(t)\xi_i\xi_j\bar{\hat{v}}+\\
-\imath\sum_{i=1}^n\bar{b}_i(t)\xi_i\bar{\hat{v}}+\bar{c}(t)\bar{\hat{v}}=0\,.
\end{multline}
Multiplying (\ref{eq_per_v_hat}) by $(t+\tau)\partial_t\bar{\hat{v}}$ and (\ref{eq_per_v_hat_con}) by $(t+\tau)\partial_t\hat{v}$ and summing the two terms yields
\begin{multline}\label{eq_sommata}
2(t+\tau)\vert\partial_t\hat{v}\vert^2+\frac{\alpha}{2}(t+\tau)\vert\xi\vert^2\omega\left(\frac{1}{\vert\xi\vert^2+1}\right)(\hat{v}\partial_t\bar{\hat{v}}+\bar{\hat{v}}\partial_t\hat{v})-\gamma(t+\tau)(\hat{v}\partial_t\bar{\hat{v}}+\bar{\hat{v}}\partial_t\hat{v})+\\
+(t+\tau)\phi_\lambda^\prime\left(\frac{t+\tau}{\beta}\right)(\hat{v}\partial_t\bar{\hat{v}}+\bar{\hat{v}}\partial_t\hat{v})-(t+\tau)\sum_{i,j=1}^na_{i,j}(t)\xi_i\xi_j(\hat{v}\partial_t\bar{\hat{v}}+\bar{\hat{v}}\partial_t\hat{v})+\\
-2(t+\tau)\sum_{i=1}^n\xi_i\Im\{b_i(t)\hat{v}\partial_t\bar{\hat{v}}\}+2(t+\tau)\Re\left\{c(t)\hat{v}\partial_t\bar{\hat{v}}\right\}=0\,.
\end{multline}
Substituting in the second term the explicit expressions of $\partial_t\hat{v}$ and $\partial_t\bar{\hat{v}}$, that may be obtained from (\ref{eq_per_v_hat}) and (\ref{eq_per_v_hat_con}), one obtains
\begin{multline}\label{eq_prima_della_grande}
2(t+\tau)\vert\partial_t\hat{v}\vert^2-\frac{\alpha^2}{2}(t+\tau)\vert\xi\vert^4\left[\omega\left(\frac{1}{\vert\xi\vert^2+1}\right)\right]^2\vert\hat{v}\vert^2+\\
+\alpha\gamma(t+\tau)\vert\xi\vert^2\omega\left(\frac{1}{\vert\xi\vert^2+1}\right)\vert\hat{v}\vert^2-\alpha(t+\tau)\vert\xi\vert^2\omega\left(\frac{1}{\vert\xi\vert^2+1}\right)\phi_\lambda^\prime\left(\frac{t+\tau}{\beta}\right)\vert\hat{v}\vert^2+\\
+\alpha(t+\tau)\vert\xi\vert^2\omega\left(\frac{1}{\vert\xi\vert^2+1}\right)\vert\hat{v}\vert^2\left(\sum_{i,j=1}^na_{i,j}(t)\xi_i\xi_j-c(t)\right)+\\
-\gamma(t+\tau)(\hat{v}\partial_t\bar{\hat{v}}+\bar{\hat{v}}\partial_t\hat{v})+(t+\tau)\phi_\lambda^\prime\left(\frac{t+\tau}{\beta}\right)(\hat{v}\partial_t\bar{\hat{v}}+\bar{\hat{v}}\partial_t\hat{v})+\\
-(t+\tau)(\hat{v}\partial_t\bar{\hat{v}}+\bar{\hat{v}}\partial_t\hat{v})\sum_{i,j=1}^na_{i,j}(t)\xi_i\xi_j+\\
-2(t+\tau)\sum_{i=1}^n\xi_i\Im\{b_i(t)\hat{v}\partial_t\bar{\hat{v}}\}+2(t+\tau)\Re\left\{c(t)\hat{v}\partial_t\bar{\hat{v}}\right\}=0\,.
\end{multline}

Integrating (\ref{eq_prima_della_grande}) between $0$ and $s$, with $s\le\sigma=1/\alpha$, yields

%Assume, now, 
%To 
%Estendiamo, ora, il valore delle funzioni $a_{i,j}$ a tutto l'asse reale, ponendo $a_{i,j}(t)=a_{i,j}(0)$ se $t<0$ e $a_{i,j}(T)$ se $t>T$ e definiamo
%$$
%a_{i,j}^{\epsilon}(t)\triangleq (\rho_{\epsilon}\ast a_{i,j})(t)=\int_{\mathbb{R}^n}\rho_{\epsilon}(t-s)a_{i,j}(s)ds\,
%$$
%dove $\rho_{\epsilon}$ \`e un mollificatore di classe $\mathscr{C}^\infty$.

\begin{multline}\label{eq_grande_prima}
2\int_0^s(t\!+\!\tau)\vert\partial_t\hat{v}(t,\xi)\vert^2dt-\frac{\alpha^2}{2}\vert\xi\vert^4\left[\omega\left(\frac{1}{\vert\xi\vert^2+1}\right)\right]^2\int_0^s(t\!+\!\tau)\vert\hat{v}(t,\xi)\vert^2dt+\\
+\underbrace{\alpha\gamma\vert\xi\vert^2\omega\left(\frac{1}{\vert\xi\vert^2+1}\right)\int_0^s(t\!+\!\tau)\vert\hat{v}(t,\xi)\vert^2dt}_{(A)}+\\
-\alpha\vert\xi\vert^2\omega\left(\frac{1}{\vert\xi\vert^2+1}\right)\int_0^s(t\!+\!\tau)\phi_\lambda^\prime\left(\frac{t\!+\!\tau}{\beta}\right)\vert\hat{v}(t,\xi)\vert^2dt+\\
+\alpha\vert\xi\vert^2\omega\left(\frac{1}{\vert\xi\vert^2+1}\right)\int_0^s(t\!+\!\tau)\sum_{i,j=1}^na_{i,j}(t)\xi_i\xi_j\vert\hat{v}(t,\xi)\vert^2dt+\\
-\alpha\vert\xi\vert^2\omega\left(\frac{1}{\vert\xi\vert^2+1}\right)\int_0^s(t\!+\!\tau)c(t)\vert\hat{v}(t,\xi)\vert^2dt+\\
+\gamma\int_0^s\vert\hat{v}(t,\xi)\vert^2dt-\gamma(s\!+\!\tau)\vert\hat{v}(s,\xi)\vert^2+\\+\underbrace{\gamma\tau\vert\hat{v}(0,\xi)\vert^2}_{(B)}+\int_0^s\left[-\phi_\lambda^{\prime\prime}\left(\frac{t\!+\!\tau}{\beta}\right)\left(\frac{t\!+\!\tau}{\beta}\right)-\phi_\lambda^\prime\left(\frac{t\!+\!\tau}{\beta}\right)\right]\vert\hat{v}(t,\xi)\vert^2dt+\\
+\underbrace{\phi_\lambda^\prime\left(\frac{s\!+\!\tau}{\beta}\right)(s\!+\!\tau)\vert\hat{v}(s,\xi)\vert^2}_{(C)}-\phi_\lambda^\prime\left(\frac{\tau}{\beta}\right)\tau\vert\hat{v}(0,\xi)\vert^2+\\
-\underbrace{\int_0^s(t\!+\!\tau)[\hat{v}(t,\xi)\partial_t\bar{\hat{v}}(t,\xi)+\bar{\hat{v}}(t,\xi)\partial_t\hat{v}(t,\xi)]\sum_{i,j=1}^na_{i,j}(t)\xi_i\xi_jdt}_{(D)}+\\
-2\sum_{i=1}^n\xi_i\int_0^s(t+\tau)\Im\{b_i(t)\hat{v}(t,\xi)\partial_t\bar{\hat{v}}(t,\xi)\}dt+\\
+2\int_0^s(t\!+\!\tau)\Re\{c(t)\hat{v}(t,\xi)\partial_t\bar{\hat{v}}(t,\xi)\}dt=0\,,
\end{multline}
where, to ease the following reasoning, some terms have been identified with capital letters from $A$ to $D$. Terms $(A)$ and $(B)$ are positive and, since $\phi$ is strictly increasing, also $(C)$ is positive. To obtain the final estimate, equation (\ref{eq_grande_prima}) needs to be slightly modified. In particular, extend functions $a_{i,j}$ to the whole real axis by setting $a_{i,j}(t)=a_{i,j}(0)$ for $t<0$ and $a_{i,j}(t)=a_{i,j}(T)$ if $t>T$ and define
$$
a_{i,j}^{\epsilon}(t)\triangleq (\rho_{\epsilon}\ast a_{i,j})(t)=\int_{\mathbb{R}^n}\rho_{\epsilon}(t-s)a_{i,j}(s)ds\,
$$
where $\rho_{\epsilon}$ is a $\mathscr{C}^\infty$ mollifier.

From (\ref{eq_grande_prima}), replacing, in $(D)$, $a_{i,j}(t)$ with $a_{i,j}(t)+a_{i,j}^\epsilon(t)-a_{i,j}^\epsilon(t)$, yields
\begin{multline}\label{eq_grande_prima_e_mezza}
2\underbrace{\int_0^s(t\!+\!\tau)\vert\partial_t\hat{v}(t,\xi)\vert^2dt}_{(E)}-\underbrace{\frac{\alpha^2}{2}\vert\xi\vert^4\left[\omega\left(\frac{1}{\vert\xi\vert^2+1}\right)\right]^2\int_0^s(t\!+\!\tau)\vert\hat{v}(t,\xi)\vert^2dt}_{(F)}+\\
-\underbrace{\alpha\vert\xi\vert^2\omega\left(\frac{1}{\vert\xi\vert^2+1}\right)\int_0^s(t\!+\!\tau)\phi_\lambda^\prime\left(\frac{t\!+\!\tau}{\beta}\right)\vert\hat{v}(t,\xi)\vert^2dt}_{(G)}+\\
+\underbrace{\alpha\vert\xi\vert^2\omega\left(\frac{1}{\vert\xi\vert^2+1}\right)\sum_{i,j=1}^n\xi_i\xi_j\int_0^s(t\!+\!\tau)a_{i,j}(t)\vert\hat{v}(t,\xi)\vert^2dt}_{(H)}+\\
-\underbrace{\alpha\vert\xi\vert^2\omega\left(\frac{1}{\vert\xi\vert^2+1}\right)\int_0^s(t\!+\!\tau)c(t)\vert\hat{v}(t,\xi)\vert^2dt}_{(I)}+\underbrace{\gamma\int_0^s\vert\hat{v}(t,\xi)\vert^2dt}_{(L)}+\\
-\underbrace{\gamma(s\!+\!\tau)\vert\hat{v}(s,\xi)\vert^2}_{(M)}+\underbrace{\int_0^s\left[-\phi_\lambda^{\prime\prime}\left(\frac{t\!+\!\tau}{\beta}\right)\!\!\left(\frac{t\!+\!\tau}{\beta}\right)\!-\!\phi_\lambda^\prime\left(\frac{t\!+\!\tau}{\beta}\right)\right]\vert\hat{v}(t,\xi)\vert^2dt}_{(N)}+\\
-\underbrace{\phi_\lambda^\prime\left(\frac{\tau}{\beta}\right)\tau\vert\hat{v}(0,\xi)\vert^2}_{(O)}+\underbrace{2\sum_{i,j=1}^n\xi_i\xi_j\int_0^s(t\!+\!\tau)\Re\{\hat{v}(t,\xi)\partial_t\bar{\hat{v}}(t,\xi)\}\widetilde{a}_{i,j}^\epsilon(t)dt}_{(P)}+\\
+\underbrace{\sum_{i,j=1}^n\xi_i\xi_j\int_0^s\vert\hat{v}(t,\xi)\vert^2\frac{\partial}{\partial t}[(t\!+\!\tau)a^\epsilon_{i,j}(t)]dt}_{(Q)}+\underbrace{\tau\sum_{i,j=1}^na^\epsilon_{i,j}(0)\xi_i\xi_j\vert \hat{v}(0,\xi)\vert^2}_{(R)}+\\
-\underbrace{(s\!+\!\tau)\sum_{i,j=1}^na^\epsilon_{i,j}(s)\xi_i\xi_j\vert \hat{v}(s,\xi)\vert^2}_{(S)}-\underbrace{2\sum_{i=1}^n\xi_i\!\int_0^s(t\!+\!\tau)\Im\{b_i(t)\hat{v}(t,\xi)\partial_t\bar{\hat{v}}(t,\xi)\}dt}_{(T)}\\
+\underbrace{2\int_0^s(t+\tau)\Re\{c(t)\hat{v}(t,\xi)\partial_t\bar{\hat{v}}(t,\xi)\}dt}_{(U)}\le 0\,,
\end{multline}
where $\widetilde{a}_{i,j}^\epsilon =a_{i,j}^\epsilon- a_{i,j}$ for all $i,j=1,\ldots,n$.

In the following each term is considered individually, beginning with $(P)$. The properties of the modulus of continuity $\omega$ guarantee that there exists a constant $C_0$ such that 
$$
\vert a_{i,j}^\epsilon(t)-a_{i,j}(t)\vert\le C_0 \omega(\epsilon)\,,
$$
for all $\epsilon$, for all $i$, for all $j$ and for all $t$. Hence
$$
\left\vert\sum_{i,j=1}^n[a_{i,j}^\epsilon(t)-a_{i,j}(t)]\xi_i\xi_j\right\vert\le \sum_{i,j=1}^n\vert a_{i,j}^\epsilon(t)-a_{i,j}(t)\vert\vert\xi_i\xi_j\vert\le C_0n^2 \omega(\epsilon)\vert\xi\vert^2\,,
$$
where the property that, for all $i$, $\vert \xi_i\vert\le\vert \xi\vert$ has been exploited. As a consequence, if
$$
\epsilon=\frac{1}{\vert\xi\vert^2+1}\,,
$$
then
$$
\vert (P)\vert\le 2C_0n^2\vert\xi\vert^2\omega\left(\frac{1}{\vert\xi\vert^2+1}\right)\int_0^s(t\!+\!\tau)\vert\hat{v}(t,\xi)\partial_t\bar{\hat{v}}(t,\xi)\vert dt\,.
$$
Young's inequality yields
$$
\vert(P)\vert\le\int_0^s(t+\tau)\vert\partial_t\hat{v}(t,\xi)\vert^2dt+C_0^2n^4\vert\xi\vert^4\left[\omega\left(\frac{1}{\vert\xi\vert^2+1}\right)\right]^2\int_0^s(t+\tau)\vert\hat{v}(t,\xi)\vert^2dt
$$
and, consequently, since $\omega(s)\in[0,1]$ for all $s\in[0,1]$ and, in turn, $-\omega(s)^2>-\omega(s)$ for all $s\in[0,1]$,
$$
(P)\ge -\underbrace{\int_0^s(t+\tau)\vert\partial_t\hat{v}(t,\xi)\vert^2dt}_{(P_1)}-\underbrace{C_0^2n^4\vert\xi\vert^4\omega\left(\frac{1}{\vert\xi\vert^2+1}\right)\int_0^s(t+\tau)\vert\hat{v}(t,\xi)\vert^2dt}_{(P_2)}\,.
$$
%Note that $(P_1)=(E)$; a first result is then
%\begin{equation}\label{eq_prima_lettere}
%2(E)+(P)\ge 2(E)-(P_1)-(P_2)=(E)-(P_2)\,.
%\end{equation}
Let us consider now the term $(Q)$. For the properties of the modulus of continuity, there exists $C_1$ such that
$$
\vert (a_{i,j}^\epsilon)^\prime(t)\vert\le C_1\frac{\omega(\epsilon)}{\epsilon}\,,
$$
for all $\epsilon$, for all $i$, for all $j$ and for all $t$. As a consequence, if
$$
\epsilon=\frac{1}{\vert\xi\vert^2+1}\,,
$$
then
\begin{multline*}
(Q)=\sum_{i,j=1}^n\xi_i\xi_j\int_0^s\vert\hat{v}(t,\xi)\vert^2(t+\tau)(a_{i,j}^\epsilon)^\prime(t)dt+\\
+\sum_{i,j=1}^n\xi_i\xi_j\int_0^s\vert\hat{v}(t,\xi)\vert^2a_{i,j}^\epsilon(t)dt\ge\\
\ge -\underbrace{C_1n^2\vert\xi\vert^2(\vert\xi\vert^2+1)\omega\left(\frac{1}{\vert\xi\vert^2+1}\right)\int_0^s(t+\tau)\vert\hat{v}(t,\xi)\vert^2dt}_{(Q_1)}+\\
+\underbrace{\sum_{i,j=1}^n\xi_i\xi_j\int_0^s a_{i,j}^\epsilon(t)\vert\hat{v}(t,\xi)\vert^2dt}_{(Q_2)}\,.
\end{multline*}
%Note that
%$$
%(H)\ge \alpha k_A \vert\xi\vert^4\omega\left(\frac{1}{\vert\xi\vert^2+1}\right)\int_0^s(t+\tau)\vert \hat{v}(t,\xi)\vert^2 dt\,,
%$$
%and
%$$
%(Q_2)\ge k_A\vert\xi\vert^2\int_0^s\vert\hat{v}(t,\xi)\vert^2dt\,.
%$$
As far as the terms (T) and (U) are concerned,
$$
(U)-(T)\ge -(U_1)-(U_2)-(T_1)-(T_2)\,,
$$
where
$$
(U_1)=2k_C^2\int_0^s(t+\tau)\vert\hat{v}(t,\xi)\vert^2dt\,,\quad (U_2)=\frac{1}{2}\int_0^s(t+\tau)\vert\partial_t\hat{v}(t,\xi)\vert^2dt\,,
$$
$$
(T_1)=2n^2k_B^2\vert\xi\vert^2\int_0^s(t+\tau)\vert\hat{v}(t,\xi)\vert^2dt\,,\quad (T_2)=\frac{1}{2}\int_0^s(t+\tau)\vert\partial_t\hat{v}(t,\xi)\vert^2dt\,.
$$
Note, moreover, that
$$
(H)\ge \alpha k_A \vert\xi\vert^4\omega\left(\frac{1}{\vert\xi\vert^2+1}\right)\int_0^s(t+\tau)\vert \hat{v}(t,\xi)\vert^2 dt\,,
$$
and
$$
(Q_2)\ge k_A\vert\xi\vert^2\int_0^s\vert\hat{v}(t,\xi)\vert^2dt\,.
$$

%As a consequence,
%\begin{multline*}
%2(E)+\frac{1}{2}(H)+\frac{1}{2}(L)+(P)+(Q)-(T)+(U)\ge\\[2mm]
% 2(E)+\frac{1}{2}(H)+\frac{1}{2}(L)-(P_1)-(P_2)-(Q_1)+(Q_2)-(T_1)-(T_2)-(U_1)-(U_2)\ge\\[2mm]
%\frac{1}{2}(H)+\frac{1}{2}(L)-(P_2)-(Q_1)+(Q_2)-(T_1)-(U_1)\,.
%\end{multline*}

\bigskip
\noindent
We claim now that there exist two positive constants $\alpha_1$ and $\gamma_1$ such that, for all $\xi\in{\mathbb R}^n$, 
\begin{multline}
\frac{\gamma_1}{4T}+\frac{\alpha_1}{2}k_A\vert\xi\vert^4\omega\left(\frac{1}{\vert\xi\vert^2+1}\right)-C_0^2n^4\vert\xi\vert^4\omega\left(\frac{1}{\vert\xi\vert^2+1}\right)+\\[2mm]
-C_1n^2\vert\xi\vert^2\left(\vert\xi\vert^2+1\right)\omega\left(\frac{1}{\vert\xi\vert^2+1}\right)-2n^2k_B^2\vert\xi\vert^2-2k_C^2+\\[2mm]
-\frac{\alpha_1^2}{2}\vert\xi\vert^4\left(\omega\left(\frac{1}{\vert\xi\vert^2+1}\right)\right)^2-\alpha_1\vert\xi\vert^2\omega\left(\frac{1}{\vert\xi\vert^2+1}\right)k_C\ge 0\,.\label{eq_per_alfa_e_gamma_prima}
\end{multline}
%and
%\begin{multline}
%\frac{\gamma^\ast}{4}+\frac{\alpha}{4}k_A\vert\xi\vert^4\omega\left(\frac{1}{\vert\xi\vert^2+1}\right)-\frac{\alpha^2}{2}\vert\xi\vert^4\left(\omega\left(\frac{1}{\vert\xi\vert^2+1}\right)\right)^2-\alpha\vert\xi\vert^2\omega\left(\frac{1}{\vert\xi\vert^2+1}\right)k_C\ge 0\,.\label{eq_per_alfa_e_gamma_seconda}
%\end{multline}
Letting the the proof of  (\ref{eq_per_alfa_e_gamma_prima}) to the reader, we 
remark that it relies on the following facts: 
when $\vert\xi\vert\ge 1$, the function
$$
\xi\to\vert\xi\vert^2\omega\left(\frac{1}{\vert\xi\vert^2+1}\right)
$$
is bounded from below by a positive quantity and 
$$
\lim_{\vert\xi\vert\to+\infty}\omega\left(\frac{1}{\vert\xi\vert^2+1}\right)=0\,.
$$
We remark also that taking a constant $\alpha\geq \alpha_1$, the inequality (\ref{eq_per_alfa_e_gamma_prima}) remains true with $\alpha$ at the place of $\alpha_1$, provided the choice of a possibly bigger $\gamma_1$.
As a consequence, if $\alpha=\max\{\alpha_1,1/T^{\prime\prime}\}$ and $\gamma\ge\gamma_1$, then
\begin{equation}\label{eq_diseg}
\frac{1}{2}(L)+\frac{1}{2}(H)-(P_2)-(Q_1)-(T_1)-(U_1)-(F)-(I)\ge 0\,.
\end{equation}

By using (\ref{eq_diseg}) %and (\ref{eq_prima_lettere}) 
into (\ref{eq_grande_prima_e_mezza}) and taking into account that $(E)=(T_2)+(U_2)=(P_1)$ and that $(R)\ge 0$, yields
\begin{multline}\label{eq_grande_terza}
\frac{1}{2}(H)+(Q_2)-\underbrace{\alpha\vert\xi\vert^2\omega\left(\frac{1}{\vert\xi\vert^2+1}\right)\int_0^s(t\!+\!\tau)\phi_\lambda^\prime\left(\frac{t+\tau}{\beta}\right)\vert\hat{v}(t,\xi)\vert^2dt}_{(G)}+\frac{1}{2}(L)+\\
-\underbrace{\gamma(s+\tau)\vert\hat{v}(s,\xi)\vert^2}_{(M)}+\underbrace{\int_0^s\left[-\phi_\lambda^{\prime\prime}\left(\frac{t\!+\!\tau}{\beta}\right)\left(\frac{t\!+\!\tau}{\beta}\right)-\phi_\lambda^\prime\left(\frac{t\!+\!\tau}{\beta}\right)\right]\vert\hat{v}(t,\xi)\vert^2dt}_{(N)}+\\
-\underbrace{\phi_\lambda^\prime\left(\frac{\tau}{\beta}\right)\tau\vert\hat{v}(0,\xi)\vert^2}_{(O)}-\underbrace{(s+\tau)\sum_{i,j=1}^na^\epsilon_{i,j}(s)\xi_i\xi_j\vert \hat{v}(s,\xi)\vert^2}_{(S)}\le 0\,.
\end{multline}

Recall, now, that $\phi_\lambda$ is a solution of equation (\ref{eq_diff_per_phi}) with $q=k_A$. Since $\omega(z)/z>1$ for all $z\in(0,1)$, equation (\ref{eq_diff_per_phi}) implies
\begin{equation}\label{seconda_eq_per_phi}
-\frac{1}{2}y\phi_\lambda^{\prime\prime}(y)>\frac{\lambda k_A}{2}\phi_\lambda^\prime(y)\,,\quad\textnormal{ for all }y\in(0,1)\,.
\end{equation}
Hence, if $\phi_\lambda$ is solution of (\ref{eq_diff_per_phi}) with $\lambda>2/k_A$,
$$
(N)\ge -\frac{1}{2}\int_0^s\phi_\lambda^{\prime\prime}\left(\frac{t+\tau}{\beta}\right)\left(\frac{t+\tau}{\beta}\right)\vert\hat{v}(t,\xi)\vert^2 dt\,,
$$
provided that $(t+\tau)/\beta\in(0,1)$ for all $t\in(0,s)$.
Consider, now, the following two cases.
\begin{enumerate}
\item If
$$
\phi_\lambda^\prime\left(\frac{t+\tau}{\beta}\right)\le \frac{(\vert\xi\vert^2+1)k_A}{4}\,,
$$
then
$$
(G)\le\frac{1}{4}\alpha k_A\vert\xi\vert^2(\vert\xi\vert^2+1)\omega\left(\frac{1}{\vert\xi\vert^2+1}\right)\int_0^s(t+\tau)\vert\hat{v}(t,\xi)\vert^2dt
$$
and hence, if
\begin{equation}\label{eq_nuova_gamma}
\gamma > \bar{\gamma}\triangleq\max\left\{\gamma_1,8T\alpha k_A\omega\left(\frac{1}{2}\right)\right\}\,,
\end{equation}
then
$$
\frac{1}{2}(H)+\frac{1}{4}(L)\ge (G)\,.
$$
In fact if $\vert\xi\vert>1$, then
$$
\frac{1}{4}\alpha k_A\vert\xi\vert^2(\vert\xi\vert^2+1)\omega\left(\frac{1}{\vert\xi\vert^2+1}\right)\int_0^s(t+\tau)\vert\hat{v}(t,\xi)\vert^2dt\le \frac{1}{2}(H)\,.
$$
If $\vert\xi\vert\le 1$, then
$$
(\vert\xi\vert^2+1)\vert\xi\vert^2\omega\left(\frac{1}{\vert\xi\vert^2+1}\right)\le 2\omega\left(\frac{1}{2}\right)
$$
and choosing $\gamma$ according to (\ref{eq_nuova_gamma}) guarantees $(G)\le (L)/4$.
\item On the contrary, if
$$
\phi_\lambda^\prime\left(\frac{t+\tau}{\beta}\right)> \frac{(\vert\xi\vert^2+1)k_A}{4}\,,
$$
then, since the function $h:(0,1)\to\mathbb{R}$ defined by $h(y)=\omega(y)/y$ is decreasing,
\begin{multline*}
(\vert\xi\vert^2+1)\omega\left(\frac{1}{\vert\xi\vert^2+1}\right)=\frac{\displaystyle \omega\left(\frac{1}{\vert\xi\vert^2+1}\right)}{\displaystyle \frac{1}{\vert\xi\vert^2+1}}\le\\[2mm]
\le\frac{\displaystyle \omega\left(\frac{k_A}{4\phi_\lambda^\prime\left(\frac{t+\tau}{\beta}\right)}\right)}{\displaystyle \frac{k_A}{4\phi_\lambda^\prime\left(\frac{t+\tau}{\beta}\right)}}=\frac{4}{k_A}\phi_\lambda^\prime\left(\frac{t+\tau}{\beta}\right)\omega\left(\frac{k_A}{4\phi_\lambda^\prime\left(\frac{t+\tau}{\beta}\right)}\right)
\end{multline*}
and, since $\omega$ is increasing,
$$
(\vert\xi\vert^2+1)\omega\left(\frac{1}{\vert\xi\vert^2+1}\right)\le\frac{4}{k_A}\phi_\lambda^\prime\left(\frac{t+\tau}{\beta}\right)\omega\left(\frac{k_A}{\phi_\lambda^\prime\left(\frac{t+\tau}{\beta}\right)}\right)\,.
$$
As a consequence, if $\phi_\lambda$ is solution of (\ref{eq_diff_per_phi}) with $\lambda>4/k_A$, then
\begin{multline}
(N)\ge-\frac{1}{2}\int_0^s\phi_\lambda^{\prime\prime}\left(\frac{t+\tau}{\beta}\right)\left(\frac{t+\tau}{\beta}\right)\vert\hat{v}(t,\xi)\vert^2dt=\\
=\frac{\lambda}{2}\int_0^s\phi_\lambda^\prime\left(\frac{t+\tau}{\beta}\right)\left(\phi_\lambda^\prime\left(\frac{t+\tau}{\beta}\right)\omega\left(\frac{k_A}{\phi_\lambda^\prime\left(\frac{t+\tau}{\beta}\right)}\right)\right)\vert\hat{v}(t,\xi)\vert^2dt\ge\\
\ge\frac{\lambda k_A}{8}(\vert\xi\vert^2+1)\omega\left(\frac{1}{\vert\xi\vert^2+1}\right)\int_0^s\phi_\lambda^\prime\left(\frac{t+\tau}{\beta}\right)\vert\hat{v}(t,\xi)\vert^2dt\,.
\end{multline}
Moreover, if
$$
\lambda>\bar{\lambda}\triangleq\max\left(\frac{4}{k_A},\frac{16T\alpha}{k_A}\right)\,,
$$
then
$$
(N)\ge\alpha(\vert\xi\vert^2+1)\omega\left(\frac{1}{\vert\xi\vert^2+1}\right)\int_0^s(t+\tau)\phi_\lambda^\prime\left(\frac{t+\tau}{\beta}\right)\vert\hat{v}(t,\xi)\vert^2dt\ge(G)\,.
$$
\end{enumerate}
In conclusion, taking into account that $(N)\ge 0$, $(H)\ge 0$, $(L)\ge 0$ and $(G)\ge 0$, leads to the inequality
\begin{equation}\label{eq_access_finale}
\frac{1}{2}(H)+\frac{1}{4}(L)+(N)-(G)\ge 0\,.
\end{equation}
Furthermore, using (\ref{eq_access_finale}) into (\ref{eq_grande_terza}) and taking into account that
$$
\frac{1}{2}(Q_2)\ge \frac{1}{2}k_A\vert\xi\vert^2\int_0^s\vert\hat{v}(t,\xi)\vert^2dt\,,
$$
yields
\begin{multline}\label{eq_finale_v}
\left(\frac{k_A\vert\xi\vert^2}{2}+\frac{\gamma}{4}\right)\int_0^s\vert\hat{v}(t,\xi)\vert^2dt\le\\
\le \phi_\lambda^\prime\left(\frac{\tau}{\beta}\right)\tau\vert\hat{v}(0,\xi)\vert^2+(s+\tau)(\gamma+k_A^{-1}\vert\xi\vert^2)\vert\hat{v}(s,\xi)\vert^2\,.
\end{multline}
Finally, substituting (\ref{eq_v_hat}) into (\ref{eq_finale_v}) yields
\begin{multline}\label{eq_u_hat_iniziale}
\frac{1}{4}\left(k_A\vert\xi\vert^2+\gamma\right)\int_0^se^{(1-\alpha t)\vert\xi\vert^2\omega\left(\frac{1}{\vert\xi\vert^2+1}\right)}e^{2\gamma t}e^{-2\beta\phi_\lambda\left(\frac{t+\tau}{\beta}\right)}\vert\hat{u}(t,\xi)\vert^2dt\le\\
\le \phi_\lambda^\prime\left(\frac{\tau}{\beta}\right)\tau e^{\vert\xi\vert^2\omega\left(\frac{1}{\vert\xi\vert^2+1}\right)}e^{-2\beta\phi_\lambda\left(\frac{\tau}{\beta}\right)}\vert\hat{u}(0,\xi)\vert^2+\\
+(s+\tau)(\gamma+k_A^{-1}\vert\xi\vert^2)e^{(1-\alpha s)\vert\xi\vert^2\omega\left(\frac{1}{\vert\xi\vert^2+1}\right)}e^{2\gamma s}e^{-2\beta\phi_\lambda\left(\frac{s+\tau}{\beta}\right)}\vert\hat{u}(s,\xi)\vert^2\,.
\end{multline}
Equation (\ref{eq_u_hat_iniziale}) holds for all $s\in(0,\sigma]$; choosing $s=\sigma$ one obtains (\ref{eq_u_hat_seconda}).$\hfill\blacksquare$

\medskip

%\begin{multline}\label{eq_u_hat_seconda}
%\frac{1}{2}\left(k_A\vert\xi\vert^2+\gamma\right)\int_0^\sigma e^{(1-\alpha t)\vert\xi\vert^2\omega\left(\frac{1}{\vert\xi\vert^2}\right)}e^{2\gamma t}e^{-2\beta\phi\left(\frac{t+\tau}{\beta}\right)}\vert\hat{u}(t,\xi)\vert^2dt\le\\
%\le \phi^\prime\left(\frac{\tau}{\beta}\right)\tau e^{\vert\xi\vert^2\omega\left(\frac{1}{\vert\xi\vert^2}\right)}e^{-2\beta\phi\left(\frac{\tau}{\beta}\right)}\vert\hat{u}(0,\xi)\vert^2+\\
%+(\sigma+\tau)(\gamma+k_A^{-1}\vert\xi\vert^2)e^{2\gamma\sigma}e^{-2\beta\phi\left(\frac{\sigma+\tau}{\beta}\right)}\vert\hat{u}(\sigma,\xi)\vert^2\,.
%\end{multline}

\subsection{An integral estimate}
Proposition~\ref{prop_con_lettere} provides a punctual estimate of the Fourier transform of $u$ which will allow us to obtain, by integration, an analogously estimate on the norm of $u$. To obtain this result the following lemma and Definition~\ref{def_OGS} are accessory.

\begin{lemma}\label{lemma_gamma}
If $u\in\mathcal{H}_0$ is solution of (\ref{eq_L}), then there exists $\bar{\gamma}$, not depending on $\xi$, such that, for all $\xi$, $e^{2\bar{\gamma} t}\vert\hat{u}(t,\xi)\vert^2$ is (weakly) increasing in $t$.$\hfill\square$
\end{lemma}
\textbf{Proof. }We want to show that there exists $\bar{\gamma}$ such that
$$
\partial_t(e^{2\bar{\gamma} t}\hat{u}(t,\xi)\bar{\hat{u}}(t,\xi))\ge 0\,.
$$
Note that
\begin{multline}\label{eq_modulo_crescente}
\partial_t(e^{2\bar{\gamma} t}\hat{u}(t,\xi)\bar{\hat{u}}(t,\xi))=2\bar{\gamma} e^{2\bar{\gamma} t}\vert\hat{u}(t,\xi)\vert^2+\\
+e^{2\bar{\gamma} t}\partial_t(\hat{u}(t,\xi))\bar{\hat{u}}(t,\xi)+e^{2\bar{\gamma} t}\hat{u}(t,\xi)\partial_t(\bar{\hat{u}}(t,\xi))\,.
\end{multline}
From (\ref{identita_trasf_four}), multiplying by $\bar{\hat{u}}(t,\xi)$ we obtain
\begin{multline*}
\bar{\hat{u}}(t,\xi)\partial_t \hat{u}(t,\xi)=\\
=\sum_{i,j=1}^na_{i,j}(t)\xi_i\xi_j\vert\hat{u}(t,\xi)\vert^2-\imath \sum_{i=1}^nb_i(t)\xi_i\vert\hat{u}(t,\xi)\vert^2+c(t)\vert\hat{u}(t,\xi)\vert^2
\end{multline*}
and also, taking in both term the complex conjugate values,
\begin{multline*}
\hat{u}(t,\xi)\partial_t \bar{\hat{u}}(t,\xi)=\\
=\sum_{i,j=1}^na_{i,j}(t)\xi_i\xi_j\vert\hat{u}(t,\xi)\vert^2+\imath \sum_{i=1}^n\bar{b}_i(t)\xi_i\vert\hat{u}(t,\xi)\vert^2+\bar{c}(t)\vert\hat{u}(t,\xi)\vert^2
\end{multline*}
and, consequently,
\begin{multline}\label{eq_per_lemma_1}
\partial_t(e^{2\bar{\gamma} t}\hat{u}(t,\xi)\bar{\hat{u}}(t,\xi))=2\bar{\gamma} e^{2\bar{\gamma} t}\vert\hat{u}(t,\xi)\vert^2+2e^{2\bar{\gamma} t}\sum_{i,j=1}^na_{i,j}(t)\xi_i\xi_j\vert\hat{u}(t,\xi)\vert^2+\\
+2e^{2\bar{\gamma} t}\sum_{i=1}^n\Im\{b_i(t)\}\xi_i\vert\hat{u}(t,\xi)\vert^2+2e^{2\bar{\gamma} t}\Re\{c(t)\}\vert\hat{u}(t,\xi)\vert^2\ge\\
2e^{2\bar{\gamma} t}\vert\hat{u}(t,\xi)\vert^2(\bar{\gamma}+k_A\vert\xi\vert^2-nk_B\vert\xi\vert-k_C)\,.
\end{multline}
Now, if $\vert \xi\vert\ge nk_B/k_A$, then $k_A\vert\xi\vert^2>nk_B\vert\xi\vert$ and hence, if $\bar{\gamma}>k_C$, we have
$$
\bar{\gamma}+k_A\vert\xi\vert^2-nk_B\vert\xi\vert-k_C\ge 0\,.
$$
On the other hand, if $\vert \xi\vert< nk_B/k_A$, then $-\vert\xi\vert>-nk_B/k_A$ and hence $-nk_B\vert\xi\vert>-n^2k_B^2/k_A$. In conclusion, the claim holds for any $\bar{\gamma}$ such that $\bar{\gamma}>2\max\{k_C,n^2k_B^2/k_A\}$. $\hfill\blacksquare$

\medskip

\medskip

Let us, now, come back to equation (\ref{eq_u_hat_seconda}). By integrating it with respect to $\xi$, the following result can be obtained. 

\begin{prop}\label{prop_norma_nuova}
Let $\sigma$ and $\tau$ be as in Proposition~\ref{prop_con_lettere}. Set $\bar{\sigma}\triangleq \sigma/8$. There exists $C>0$ such that, whenever $u\in\mathcal{H}_0$ is a solution of (\ref{eq_L}), with $\mathcal{L}$ fulfilling Assumption~\ref{ipot_princ}, one has, for all $\beta\ge \sigma+\tau$,
\begin{multline}\label{eq_u_settima}
\sup_{z\in[0,\bar{\sigma}]}\Vert u(z,\cdot)\Vert^2_{H_{\frac{1}{2},\omega}^1}\le\\[2mm]
\le C e^{-\sigma\phi^\prime\left(\frac{\sigma+\tau}{\beta}\right)}\left[\phi^\prime\left(\frac{\tau}{\beta}\right)e^{-2\beta\phi\left(\frac{\tau}{\beta}\right)}\Vert u(0,\cdot)\Vert^2_{H_{1,\omega}^0}+\Vert u(\sigma,\cdot)\Vert^2_{H^1}\right]\,,
\end{multline}
\end{prop}
where $\phi=\phi_{\bar{\lambda},k_A}$ with $\bar{\lambda}$ given by Proposition~\ref{prop_con_lettere}. The constant $C$ depends no $n$, $k_A$, $k_B$, $k_C$, $\omega$, $T$ and $T^{\prime\prime}$.
$\hfill\square$

\textbf{Proof. }In the hypotheses of the claim, Proposition~\ref{prop_con_lettere} guarantees the existence of $\sigma$, $\alpha$, $\gamma$ and $\phi_{\lambda}$ such that (\ref{eq_u_hat_seconda}) holds. The integrand function in (\ref{eq_u_hat_seconda}) is positive and, consequently, the term on the left hand side can be bounded from below by integrating on an interval contained in $[0,\sigma]$. Let $\tau\le \sigma/4$ and let $z$ be a value such that $0<z\le \bar{\sigma}$; we have
$$
[z,2z+\tau]\subset [0,\sigma/2]\,;
$$
by integrating with respect to $\xi$ and taking into account that, since $\sigma=1/\alpha$,
$$
1-\alpha t\ge 1-\alpha \frac{\sigma}{2}\ge\frac{1}{2}\,,
$$
for all $t\in[0,\sigma/2]$, one obtains
\begin{multline}\label{eq_u_hat_terza}
\frac{1}{4}\int_{\mathbb{R}^n}\left(k_A\vert\xi\vert^2+\gamma\right)e^{\frac{1}{2}\vert\xi\vert^2\omega\left(\frac{1}{\vert\xi\vert^2+1}\right)}\int_z^{2z+\tau}e^{2\gamma t}e^{-2\beta\phi_{\lambda}\left(\frac{t+\tau}{\beta}\right)}\vert\hat{u}(t,\xi)\vert^2dtd\xi\le\\
\le \tau \phi_{\lambda}^\prime\left(\frac{\tau}{\beta}\right)e^{-2\beta\phi_{\lambda}\left(\frac{\tau}{\beta}\right)}\int_{\mathbb{R}^n} e^{\vert\xi\vert^2\omega\left(\frac{1}{\vert\xi\vert^2+1}\right)}\vert\hat{u}(0,\xi)\vert^2d\xi+\\
+(\sigma+\tau)e^{2\gamma\sigma}e^{-2\beta\phi_{\lambda}\left(\frac{\sigma+\tau}{\beta}\right)}\int_{\mathbb{R}^n}(\gamma+k_A^{-1}\vert\xi\vert^2)\vert\hat{u}(\sigma,\xi)\vert^2d\xi\,.
\end{multline}

\medskip

Now, let $\bar{\bar{\gamma}}$ be a value of $\gamma$ fulfilling equation (\ref{eq_nuova_gamma}), let $\bar{\gamma}$ be the value provided by Lemma~\ref{lemma_gamma} and let
$$
\gamma>\max\{\bar{\bar{\gamma}},\bar{\gamma}\}\,.
$$
Since $\phi_{\lambda}$ is increasing, we have that
$$
e^{-2\beta\phi_{\lambda}\left(\frac{t+\tau}{\beta}\right)}\ge e^{-2\beta\phi_{\lambda}\left(\frac{2(z+\tau)}{\beta}\right)}
$$
for all $t<2z+\tau$. As a consequence, using also the fact that $e^{2\gamma z}\ge 1$, equation (\ref{eq_u_hat_terza}) yields
\begin{multline}\label{eq_u_hat_quarta}
c_1(z+\tau)\int_{\mathbb{R}^n}(\vert\xi\vert^2+1)e^{\frac{1}{2}\vert\xi\vert^2\omega\left(\frac{1}{\vert\xi\vert^2+1}\right)}\vert\hat{u}(z,\xi)\vert^2d\xi\le\\
\le \tau\phi_{\lambda}^\prime\left(\frac{\tau}{\beta}\right)e^{2\beta\left[\phi_{\lambda}\left(\frac{2(z+\tau)}{\beta}\right)-\phi_{\lambda}\left(\frac{\tau}{\beta}\right)\right]}\int_{\mathbb{R}^n} e^{\vert\xi\vert^2\omega\left(\frac{1}{\vert\xi\vert^2+1}\right)}\vert\hat{u}(0,\xi)\vert^2d\xi+\\
+c_2(\sigma+\tau)e^{2\gamma\sigma}e^{2\beta\left[\phi_{\lambda}\left(\frac{2(z+\tau)}{\beta}\right)-\phi_{\lambda}\left(\frac{\sigma+\tau}{\beta}\right)\right]}\int_{\mathbb{R}^n}(1+\vert\xi\vert^2)\vert\hat{u}(\sigma,\xi)\vert^2d\xi\,,
\end{multline}
where the constant values
$$
c_1\triangleq\frac{1}{4}\min\left\{k_A,\gamma\right\}\,,\qquad c_2\triangleq \max\left\{\gamma,k_A^{-1}\right\}
$$
have been introduced.
%Definiamo ora, per una funzione $u\in\mathscr{C}^\omega[0,T]\cup L^\infty([0,T])$ la norma
%\begin{equation}\label{eq_norma}
%\Vert u \Vert_{H_{q,\omega}^\nu}=\left(\int_{\mathbb{R}^n}(1+\vert\xi\vert^2)^\nu e^{q\vert\xi\vert^2\omega\left(\frac{1}{\vert\xi\vert^2}\right)}\vert \hat{u}(\xi)\vert^2d\xi\right)^2\,.
%\end{equation}
%Con l'uso della norma (\ref{eq_norma}), 
Dividing by $\tau$ and taking into account that $(z+\tau)/\tau>1$ and that $\phi_{\lambda}$ is negative, it is easy to see that (\ref{eq_u_hat_quarta}) implies
\begin{multline}\label{eq_u_quinta}
c_1\Vert u(z,\cdot)\Vert^2_{H_{\frac{1}{2},\omega}^1}\le \phi_{\lambda}^\prime\left(\frac{\tau}{\beta}\right)e^{2\beta\left[\phi_{\lambda}\left(\frac{2(z+\tau)}{\beta}\right)-\phi_{\lambda}\left(\frac{\tau}{\beta}\right)\right]}\Vert u(0,\cdot)\Vert^2_{H_{1,\omega}^0}+\\
+c_2\frac{\sigma+\tau}{\tau}e^{2\gamma\sigma}e^{2\beta\left[\phi_{\lambda}\left(\frac{2(z+\tau)}{\beta}\right)-\phi_{\lambda}\left(\frac{\sigma+\tau}{\beta}\right)\right]}\Vert u(\sigma,\cdot)\Vert^2_{H^1}\le\\
\le\phi_{\lambda}^\prime\left(\frac{\tau}{\beta}\right)e^{2\beta\left[\phi_{\lambda}\left(\frac{2(z+\tau)}{\beta}\right)-\phi_{\lambda}\left(\frac{\tau}{\beta}\right)-\phi_{\lambda}\left(\frac{\sigma+\tau}{\beta}\right)\right]}\Vert u(0,\cdot)\Vert^2_{H_{1,\omega}^0}+\\
+c_2\frac{\sigma+\tau}{\tau}e^{2\gamma\sigma}e^{2\beta\left[\phi_{\lambda}\left(\frac{2(z+\tau)}{\beta}\right)-\phi_{\lambda}\left(\frac{\sigma+\tau}{\beta}\right)\right]}\Vert u(\sigma,\cdot)\Vert^2_{H^1}\,,
\end{multline}
Moreover, with respect to $\phi_{\lambda}$, note that since $\phi_{\lambda}$ is increasing,
$$
2(z+\tau)\le\frac{\sigma}{2}+\tau\quad\Rightarrow\quad\phi_{\lambda}\left(\frac{2(z+\tau)}{\beta}\right)\le\phi_{\lambda}\left(\frac{\frac{\sigma}{2}+\tau}{\beta}\right)\,.
$$
In addition, since $\phi_{\lambda}$ is also concave,
$$
\phi_{\lambda}\left(\frac{\sigma+\tau}{\beta}\right)-\phi_{\lambda}\left(\frac{\frac{\sigma}{2}+\tau}{\beta}\right)\ge\frac{\sigma}{2\beta}\phi_{\lambda}^\prime\left(\frac{\sigma+\tau}{\beta}\right)\,.
$$
As a consequence, from (\ref{eq_u_quinta}) one obtains
\begin{multline}\label{eq_u_sesta}
c_1\Vert u(z,\cdot)\Vert^2_{H_{\frac{1}{2},\omega}^1}\le\\
\le e^{-\sigma\phi_{\lambda}^\prime\left(\frac{\sigma+\tau}{\beta}\right)}\left[\phi_{\lambda}^\prime\left(\frac{\tau}{\beta}\right)e^{-2\beta\phi_{\lambda}\left(\frac{\tau}{\beta}\right)}\Vert u(0,\cdot)\Vert^2_{H_{1,\omega}^0}+c_2\frac{\sigma+\tau}{\tau}e^{2\gamma\sigma}\Vert u(\sigma,\cdot)\Vert^2_{H^1}\right]\,,
\end{multline}
namely
\begin{multline}\label{eq_u_settima_seconda}
\Vert u(z,\cdot)\Vert^2_{H_{\frac{1}{2},\omega}^1}\le\\
\le C e^{-\sigma\phi_{\lambda}^\prime\left(\frac{\sigma+\tau}{\beta}\right)}\left[\phi_{\lambda}^\prime\left(\frac{\tau}{\beta}\right)e^{-2\beta\phi_{\lambda}\left(\frac{\tau}{\beta}\right)}\Vert u(0,\cdot)\Vert^2_{H_{1,\omega}^0}+\Vert u(\sigma,\cdot)\Vert^2_{H^1}\right]\,,
\end{multline}
where
$$
C=\max\left\{\frac{1}{c_1},\frac{c_2(\sigma+\tau)e^{2\gamma\sigma}}{c_1\tau}\right\}\,.
$$
Equation (\ref{eq_u_settima_seconda}) holds for all $z\in[0,\bar{\sigma}]$ and hence equation (\ref{eq_u_settima}) immediately follows.$\hfill\blacksquare$

\medskip

\subsection{Proof of Theorem~\ref{teo_finale_stab_norme_classiche}}

Proposition~\ref{prop_norma_nuova} states, in particular, that the norm of $u$ in any insatant of the sub-interval $[0,\bar{\sigma}]\subset[0,\sigma]$ is bounded by a quantity depending on the value of the norm in the initial and final instants, i.e. on $\Vert u(0,\cdot)\Vert_{H_{1,\omega}^0}$ and $\Vert u(\sigma,\cdot)\Vert_{H^1}$. Nevertheless, to obtain a stability result, the right hand side term in equation (\ref{eq_u_settima_seconda}) must tend to zero when $\Vert u(0,\cdot)\Vert_{H_{1,\omega}^0}$ tends to zero, which is not immediate to guess. The following lemma allows one to choose $\beta$ in such a way that (\ref{eq_u_settima_seconda}) can be written in a form from which the stability property can be obtained more easily.

%
%Supponiamo ora -- lo verificheremo tra breve -- che si possa sempre scegliere $\beta$ in modo da soddisfare l'equazione
%\begin{equation}\label{eq_beta}
%e^{-2\beta\phi\left(\frac{\tau}{\beta}\right)}\phi^\prime\left(\frac{\tau}{\beta}\right)=\Vert u(0)\Vert_{H^0_{1,\omega}}^{-2}\,.
%\end{equation}
%Con questa scelta di $\beta$, la (\ref{eq_u_settima}) diventa
%\begin{equation}\label{eq_u_ottava}
%\Vert u(z)\Vert^2_{H_{\frac{1}{2},\omega}^1}\le c_3e^{-\sigma\phi^\prime\left(\frac{\sigma+\tau}{\beta}\right)}[1+\Vert u(\sigma)\Vert^2_{H^1}]\,.
%\end{equation}
%Per concludere con il risultato di stabilit\`a, ci serviremo dei seguenti due lemmi.
\begin{lemma}\label{lem_finale}
Let $\phi$ be a solution of (\ref{eq_diff_per_phi}) with $\lambda>0$ and $q>0$ and let $\tau>0$. Let $h:]0,1[\to]q,+\infty[$ be defined by
$$
h(z)\triangleq e^{-2\tau \phi(z)/z}\phi^\prime(z)\,.
$$
The function $h$ so defined is strictly decreasing with
$$
\lim_{z\to 0}h(z)=+\infty\,,\quad \lim_{z\to 1}h(z)=q\,.
$$
$\hfill\square$
%and for all $y<1/q$, there exists a unique $\beta>\tau$ such that
%\begin{equation}\label{eq_beta}
%e^{-2\beta\phi\left(\frac{\tau}{\beta}\right)}\phi^\prime\left(\frac{\tau}{\beta}\right)=\frac{1}{y}\,.
%\end{equation}
%Moreover, denoting by $\beta^\ast$ the function which associates with $y$ the maximum between the solution to (\ref{eq_beta}) and the value $T+\tau$, we have
%\begin{equation}\label{eq_beta_sconda_prop}
%\lim_{y\to 0}e^{-\sigma\phi^\prime\left(\frac{\sigma+\tau}{\beta^\ast(y)}\right)}=0\,.
%\end{equation}
%for all $\sigma\in(0,T]$.$\hfill\square$
\end{lemma}
\textbf{Proof. }The claim is easily proven by computing $h^\prime$.$\hfill\blacksquare$
%\textbf{Proof. }Equation (\ref{eq_beta}) is equivalent to
%$$
%2\beta\phi\left(\frac{\tau}{\beta}\right)-\log\left(\phi^\prime\left(\frac{\tau}{\beta}\right)\right)=\log(y)\,.
%$$
%Let $f_\tau:[\tau,+\infty]\to\mathbb{R}$ be the function defined by
%$$
%f_\tau(\beta)=2\beta\phi\left(\frac{\tau}{\beta}\right)-\log\left(\phi^\prime\left(\frac{\tau}{\beta}\right)\right)\,.
%$$
%It is easy to obtain (see Section~\ref{sec_omega})
%$$
%\lim_{\beta\to+\infty}f_\tau(\beta)=-\infty\,,\qquad\lim_{\beta\to\tau^+}f_\tau(\beta)=-\log(\phi^\prime(1))=-\log(q)\,.
%$$
%Moreover,
%$$
%\frac{\partial f_\tau(\beta)}{\partial \beta}=2\phi\left(\frac{\tau}{\beta}\right)-\frac{2\tau}{\beta}\phi^\prime\left(\frac{\tau}{\beta}\right)+\frac{\tau}{\beta^2\phi^\prime\left(\frac{\tau}{\beta}\right)}\phi^{\prime\prime}\left(\frac{\tau}{\beta}\right)<0\,.
%$$
%As a consequence, $f_\tau$ is monotonically decreasing and takes on all the values between $-\infty$ and $-\log q$. Now, the assumption on $y$ guarantees that $\log y<-\log q$; hence there always exists $\beta>\tau$ such that $f_\tau(\beta)=\log y$. The first statement is so proven. Note, now, that the function $\widehat{\beta}:(0,1/q)\to(\tau,+\infty)$ defined by $\widehat{\beta}(y)=f_\tau^{-1}(\log y)$, is such that
%$$
%\lim_{y\to 0}\beta^\ast(y)=\lim_{y\to 0}\widehat{\beta}(y)=+\infty\,.
%$$
%To conclude the proof it is sufficient to note that
%$$
%\lim_{\beta\to+\infty}e^{-\sigma\phi^\prime\left(\frac{\sigma+\tau}{\beta}\right)}=0\,,
%$$
%for all $\tau>0$ and for all $\sigma>0$.$\hfill\blacksquare$

As a consequence of Lemma~\ref{lem_finale}, $h$ can be inverted and its inverse $h^{-1}:]q,+\infty[\to ]0,1[$ is strictly increasing and
$$
\lim_{y\to+\infty}h^{-1}(y)=0\,.
$$
\medskip

Now the main stability result can be proven.

\medskip

\textbf{Proof of Theorem~\ref{teo_finale_stab_norme_classiche}.}
In (\ref{eq_u_settima}) of Proposition~\ref{prop_norma_nuova} we want to choose $\beta>\sigma+\tau$ in such a way that
$$
\phi^\prime\left(\frac{\tau}{\beta}\right)e^{-2\beta\phi\left(\frac{\tau}{\beta}\right)}=\Vert u(0,\cdot)\Vert_{H_{1,\omega}^0}^{-2}\,.
$$
This goal is achieved by taking
$$
\beta=\frac{\tau}{h^{-1}\left(\Vert u(0,\cdot)\Vert^2_{H_{1,\omega}^0}\right)}\,,
$$
provided that $\Vert u(0,\cdot)\Vert_{H_{1,\omega}^0}<q^{-1/2}$ and $\Vert u(0,\cdot)\Vert_{H^0_{1,\omega}}<h\left(\frac{\tau}{\sigma+\tau}\right)^{-1/2}$. With this choice of $\beta$, one obtains, from (\ref{eq_u_settima}),
%
%
%
%Applying Lemma~\ref{lem_finale} to equation (\ref{eq_u_settima}) with $y=\Vert u(0,\cdot)\Vert^2_{H_{1,\omega}^0}$ allows one to claim that there exist $\sigma\in(0,T)$ and $C$, depending on $k_A$, $k_B$, $k_C$, $\omega$, $T$, $n$ and on the value of the $\mathscr{C}^\omega$-norm of the coefficients of the principal part of $\mathscr{L}$, and there exists $\bar{\sigma}<\sigma$ such that, if $\Vert u(0,\cdot)\Vert^2_{H_{1,\omega}^0}<1/k_A$, then \textcolor{red}{*** anche qui sigma bar e sigma/8. Inoltre, perche' possiamo prendere 1 come esponente al secondo membro? ***}
\begin{equation}\label{eq_u_settima_ultima}
\sup_{z\in[0,\bar{\sigma}]}\Vert u(z,\cdot)\Vert_{H_{\frac{1}{2},\omega}^1}^2\le C e^{-\sigma \widehat{g}\left(\Vert u(0,\cdot)\Vert^2_{H_{1,\omega}^0}\right)}\left[1+\Vert u(\sigma,\cdot)\Vert^2_{H^1}\right]\,,
\end{equation}
where $\widehat{g}$ is defined by
$$
\widehat{g}(y)=\phi^\prime\left(\frac{\sigma+\tau}{\tau}h^{-1}(y^{-1})\right)\,,
$$
so that
$$
\lim_{y\to 0}\widehat{g}(y)=+\infty\,.
$$
Note, in particular, that taking $\tau=\sigma/4$ the condition $\Vert u(0,\cdot)\Vert_{H^0_{1,\omega}}<h(\tau/(\sigma+\tau))^{-1/2}$ yields $\Vert u(0,\cdot)\Vert_{H^0_{1,\omega}}<\widehat{\rho}$ where
$$
\widehat{\rho}\triangleq\min\{e^{-\tau\frac{4}{5}\phi\left(\frac{5}{4}\right)}\phi^\prime\left(\frac{5}{4}\right)^{1/2},q^{-1/2}\}\,.
$$
Note, now, that
\begin{equation}\label{eq_ultime_1}
\Vert u(z,\cdot)\Vert_{H^1}^2\le \Vert u(z,\cdot)\Vert_{H_{\frac{1}{2},\omega}^1}^2
\end{equation}
and that, for all $\nu>0$ and all $\epsilon>0$, there exists $\widetilde{C}_{\nu,\epsilon}$ such that
$$
\Vert u(0,\cdot)\Vert_{H_{1,\omega}^0}^2\le\widetilde{C}_{\nu,\epsilon}\Vert u(0,\cdot)\Vert_{H_{\nu,\epsilon}^0}^2\,.
$$
It follows that
\begin{equation}\label{eq_ultima_questa}
\sup_{z\in[0,\bar{\sigma}]}\Vert u(z,\cdot)\Vert^2_{H_1}\le C e^{-\sigma \widehat{g}\left(\widetilde{C}_{\nu,\epsilon}\Vert u(0,\cdot)\Vert^2_{H^0_{\nu,\epsilon}}\right)}\left[1+\Vert u(\sigma,\cdot)\Vert^2_{H^1}\right]\,,
\end{equation}
provided that
$$
\Vert u(0,\cdot)\Vert_{H^0_{\nu,1}}<\frac{\widehat{\rho}}{C^{1/2}_{\nu,\epsilon}}\,.
$$

%
%Being $\widetilde{g}$ strictly decreasing, immediately yields
%\begin{equation}\label{eq_ultime_2}
%-\widetilde{g}\left(\Vert u(0,\cdot)\Vert_{H_{1,\omega}^0}\right)\le-\widetilde{g}\left(\widetilde{C}\Vert u(0,\cdot)\Vert_{H_{\nu,\epsilon}^0}\right)\,.
%\end{equation}
%In addition,
%\begin{equation}\label{eq_ultime_3}
%\Vert u(0,\cdot)\Vert^2_{H^0_{\nu,\epsilon}}\le \frac{1}{\widetilde{C}k_A}\quad\Rightarrow\quad \Vert u(0,\cdot)\Vert^2_{H^0_{1,\omega}}\le \frac{1}{k_A}\,.
%\end{equation}
By defining $g(y)=\widetilde{g}(\widetilde{C}_{\nu,\epsilon}y)$, equation (\ref{eq_ultima_questa}) allows one to easily obtain (\ref{eq_u_settima_ultima_norme_classiche}).$\hfill\blacksquare$

%Let $\mathcal{L}$ be a backward parabolic operator whose coefficients depend only on $t$, i.e. let
%$$
%\mathcal{L}u=\partial_t u+\sum_{i,j=1}^na_{i,j}(t)\partial_{x_i}\partial_{x_j}u+\sum_{j=1}^nb_j(t)\partial_{x_j}u+c(t)u
%$$
%on the strip $[0,T]\times\mathbb{R}^n$. Suppose that $a_{i,j}(t)=a_{j,i}(t)$ for all $i,j=1,\ldots,n$ and for all $t\in[0,T]$. Let $a_{i,j},b_j,c\in L^\infty([0,T])$, for all $i,j=1,\ldots,n$. Let $\mu$ be a modulus of continuity satisfying the Osgood condition. Let $a_{i,j}$ $\mathscr{C}^\mu$-continuous, i.e.
%$$
%\sup_{0<\vert\tau\vert<1}\frac{\vert a_{i,j}(t+\tau)-a_{i,j}(t)\vert}{\mu(\vert\tau\vert)}<+\infty\,.
%$$

\medskip

The claim of Theorem~\ref{teo_finale_stab_norme_classiche} to the whole interval $[0,T]$. 
%To obtain this result, Theorem~\ref{} can be iteratively applied to adjacent intervals according to the following idea. Equation () provides an estimate of the $H^1$-norm of $u$ in the interval $[0,\bar{\sigma}]$; this estimate depends, in particular, on the $H_{\nu,\epsilon}^0$-norm of $u$ in the initial instant, i.e. in $t=0$. If one could obtain an estimate in which .... NON SO SE QUESTO DISCORSO "propedeutico" SERVE

\subsection{Proof of Theorem~\ref{teo_nuovo}}
Theorem~\ref{teo_nuovo} is proven iterating a finite number of times the estimate given by the following lemma.
\begin{lemma}\label{lemma_iterativo}
Under the same hypotheses of Theorem~\ref{teo_finale_stab_norme_classiche},
$$
\sup_{z\in[0,\bar{\sigma}/2]}\Vert u(z,\cdot)\Vert_{L^2}\le C^\prime Ce^{-\sigma g\left(C^{\prime\prime}\Vert u(0,\cdot)\Vert^2_{L^2}\right)}\left[1+\Vert u(\sigma,\cdot)\Vert^2_{L^2}\right]\,.
$$
The constants $C^\prime$ and $C^{\prime\prime}$ depend on $n$, $k_A$, $k_B$, $k_C$, $\nu$, $\epsilon$ and $\sigma$ and tend to $+\infty$ as $\sigma$ tends to zero.$\hfill\square$
\end{lemma}

\textbf{Proof. }Analogously to Lemma~\ref{lem_tutti gevrey}, extend $a_{i,j}$, $b_i$ and $c$ on $[-\sigma/2,T]$ and $u$ to a solution of $\mathcal{L}$ on $[-\sigma/2,T]$. Then the results of Theorem~\ref{teo_finale_stab_norme_classiche} on $[-\bar{\sigma}/2,T-\bar{\sigma}/2]$ gives
$$
\sup_{z\in[-\bar{\sigma}/2,\bar{\sigma}/2]}\Vert u(z,\cdot)\Vert^2_{H_1}\le C e^{-\sigma g\left(\Vert u(-\bar{\sigma}/2,\cdot)\Vert^2_{H^0_{\nu,\epsilon}}\right)}\left[1+\Vert u(\sigma-\bar{\sigma}/2)\Vert^2_{H_1}\right]\,.
$$
By Lemma~\ref{lem_tutti gevrey} we obtain
\begin{multline}
\sup_{z\in[0,\bar{\sigma}/2]}\Vert u(z,\cdot)\Vert^2_{L^2}\le C e^{-\sigma g\left(C^{\prime\prime}\Vert u(0,\cdot)\Vert^2_{L^2}\right)}\left[1+\Vert u(\sigma-\frac{\sigma}{16},\cdot)\Vert^2_{H^1}\right]\le\\[2mm]
\le C^\prime C e^{-\sigma g\left(C^{\prime\prime}\Vert u(0,\cdot)\Vert^2_{L^2}\right)}\left[1+\Vert u(\sigma,\cdot)\Vert^2_{L^2}\right]\,.
\end{multline}
$\hfill\blacksquare$

Now set $G(y)\triangleq (1+D)C^\prime C e^{-\sigma g\left(C^{\prime\prime} y\right)}$ and note that $\lim_{y\to 0}G(y)=0$. We have just proven that
\begin{equation}\label{eq_con_la_G}
\sup_{z\in[0,\bar{\sigma}/2]}\Vert u(z,\cdot)\Vert^2_{L^2}\le G\left(\Vert u(0,\cdot)\Vert^2_{L^2}\right)\,.
\end{equation}
Finally, let $T^\prime:0<T^\prime<T$; take $T^{\prime\prime}=(T+T^\prime)/2$ (so that $T^\prime<T^{\prime\prime}<T$). Note that $\bar{\sigma}/2=\sigma/16$ and recall that $\sigma=\min\{1/\alpha_1,T^{\prime\prime}\}$. To complete the proof of Theorem~\ref{teo_nuovo} it is sufficient to iterate inequality (\ref{eq_con_la_G}) a finite number of times. Indeed, set $T_0=0$ and, for $i\ge 0$,
$$
T_{i+1}=T_i+\frac{1}{16}\min\left\{\frac{1}{\alpha_1},T^{\prime\prime}-T_i\right\}\,.
$$
For all $i$ inequality (\ref{eq_con_la_G}) provides
$$
\sup_{z\in[T_i,T_{i+1}]}\Vert u(z,\cdot)\Vert^2_{L^2}\le G_i\left(\Vert u(T_i,\cdot)\Vert^2_{L^2}\right)\,.
$$
The result follows by noting that
$$
T_{i+1}-T_i=\frac{1}{16}\min\left\{\frac{1}{\alpha_1},T^{\prime\prime}-T_i\right\}\,,
$$
and that, for all $j$
$$
T_{j+1}=\sum_{i=0}^j\frac{1}{16}\min\left\{\frac{1}{\alpha_1},T^{\prime\prime}-T_i\right\}\,.
$$
The sequence $\left\{T_j\right\}_{j\in\mathbb{N}}$ is increasing and bounded from above by $T^{\prime\prime}$; hence it admits a limit. Let this limit be $T^\ast$; we want to show that $T^\ast=T^{\prime\prime}$. Obviously, $T^\ast\le T^{\prime\prime}$; suppose that $T^\ast<T^{\prime\prime}$, then $T^{\prime\prime}-T_i\ge T^{\prime\prime}-T^\ast>0$ and, consequently,
$$
T_{j+1}\ge \sum_{i=0}^j\frac{1}{16}\min\left\{\frac{1}{\alpha_1},T^{\prime\prime}-T^\ast\right\}
$$
for all $j$, yielding $\lim_{j\to\infty}T_j=+\infty$, which is a contradiction. Therefore it must be $T^\ast=T^{\prime\prime}$ which means that $T_j>T^\prime$ for some $j$.$\hfill\blacksquare$

\section{A specific case}\label{sec_part}
In this section the explicit expression of the function $G$ appearing in the statement of~\ref{teo_nuovo} is computed when the modulus of continuity $\omega:]0,e^{1-e}]\to\mathbb{R}$ is defined by
$$
\omega(s)=s(1-\log s)\log(1-\log s)\,.
$$
Note that $\omega$ is increasing, fulfils the Osgood condition but is not a Log-Lipschitz function. Consider, now, the function $\theta:[e^{e-1},+\infty[\to [0,+\infty[$ defined by
\begin{equation*}%\label{eq_def_theta}
\theta(\tau)=\int_{1/\tau}^{e^{1-e}}\frac{1}{\omega(s)}ds=\log(\log(1+\log\tau))
\end{equation*}
and the function $\psi_{\lambda,q}:]0,1]\to [e^{e-1},+\infty[$ defined by
\begin{equation}\label{eq_def_psi}
\psi_{\lambda,q}(y)=\theta^{-1}(-\lambda q\log y)=\exp(e^{y^{-\lambda q}}-1)\,.
\end{equation}
From the definition of $\psi_{\lambda,q}$, one can easily check that it is strictly decreasing and that
\begin{equation}
\psi_{\lambda,q}^\prime(y)=\exp\left(e^{y^{-\lambda q}}-1\right)e^{y^{-\lambda q}}(-\lambda q)y^{-\lambda q-1}=-\frac{\lambda q}{y}\left(\psi_{\lambda,q}(y)\right)^2\omega\left(\frac{1}{\psi_{\lambda,q}(y)}\right)\,,
\end{equation}
hence the function $\phi_{\lambda,q}:]0,1]\to]-\infty,0]$ defined by
\begin{equation*}\label{eq_phi_lambda_integrale}
\phi_{\lambda,q}(y)=-q\int_y^1 \psi_{\lambda,q}(z)dz
\end{equation*}
is such that
$$
\phi_{\lambda,q}^{\prime\prime}(y)=-\frac{\lambda}{y}\left(\phi^\prime_{\lambda,q}(y)\right)^2\omega\left(\frac{q}{\phi_{\lambda,q}^\prime}(y)\right)
$$
i.e. $\phi_{\lambda,q}$ is a solution of equation (\ref{eq_diff_per_phi}). Note, as an accessory result, that 
$$
\phi_{\lambda,q}^\prime(y)=q\phi_{\lambda,q}(y)\ge q e^{e-1}\,.
$$
From now on, we choose $q=k_A$ and $\lambda\ge\bar{\lambda}$ as in the proof of Proposition~\ref{prop_con_lettere} and, for the sake of a simpler notation, we write $\phi_\lambda$ and $\psi_\lambda$ instead of $\phi_{\lambda,q}$ and $\psi_{\lambda,q}$, respectively. Proposition~\ref{prop_norma_nuova} then, gives
\begin{equation}\label{eq_stima_seconda}
\sup_{z\in[0,\bar{\sigma}]}\Vert u(z,\cdot)\Vert_{L^2}^2\le C e^{-\sigma\phi_\lambda^\prime\left(\frac{\sigma+\tau}{\beta}\right)}\phi^\prime_\lambda\left(\frac{\tau}{\beta}\right)\left[e^{-2\beta\phi_\lambda\left(\frac{\tau}{\beta}\right)}\Vert u(0,\cdot)\Vert^2_{H_{1,\omega}^0}+\Vert u(\sigma,\cdot)\Vert_{H^1}^2\right]\,.
\end{equation}
Arguing as in Lemma~\ref{lemma_iterativo} one may obtain
\begin{equation}\label{eq_sessanta}
\sup_{z\in[0,\bar{\sigma}/2]}\Vert u(z,\cdot)\Vert^2_{L^2}\le C e^{-\sigma\phi^\prime_{\lambda}\left(\frac{\sigma+\tau}{\beta}\right)}\phi_{\lambda}^\prime\left(\frac{\tau}{\beta}\right)\left[e^{-2\beta\phi_\lambda\left(\frac{\tau}{\beta}\right)}\Vert u(0,\cdot)\Vert^2_{L^2}+\Vert u(\sigma,\cdot)\Vert^2_{L^2}\right]\,.
\end{equation}
We, now, introduce the function $\Lambda:[0,+\infty[\to]-\infty,0]$ defined by
\begin{equation}\label{eq_def_Lambda}
\Lambda(y)=y\phi_\lambda\left(\frac{1}{y}\right)
\end{equation}
which is strictly decreasing and, hence, invertible. Its inverse, $\Lambda^{-1}:]-\infty,0]\to[1,+\infty[$ is also strictly decreasing. We want to find a value of $\beta>\sigma+\tau$ such that
$$
e^{-2\tau\Lambda\left(\frac{\beta}{\tau}\right)}=\Vert u(0,\cdot)\Vert^{-2}_{L^2}\,.
$$
Easy computations yield
\begin{equation}\label{eq_espressione_beta}
\beta=\tau\Lambda^{-1}\left(\frac{1}{\tau}\log\Vert u(0,\cdot)\Vert_{L^2}\right)
\end{equation}
Note that this value of $\beta$ is larger than $\sigma+\tau$ if and only if
$$
\Vert u(0\cdot)\Vert_{L^2}<e^{\tau\Lambda\left(\frac{\sigma+\tau}{\tau}\right)}\triangleq \rho\,.
$$
In particular, if $\tau=\sigma/4$ then $\rho=e^{\tau\Lambda\left(5/4\right)}$; we show below that a smaller value of $\tau$ performs better. Note, now, that for $\zeta>1$ and $y<1/\zeta$
$$
\log\left(\psi_{\lambda,q}(\zeta y)\right)=\left(\log\left(\psi_{\lambda,q}(y)\right)+1\right)^{\zeta^{-\lambda q}}-1\,;
$$
therefore
\begin{equation}\label{eq_sessanta4}
\phi^\prime_{\lambda}\left(\frac{\sigma+\tau}{\beta}\right)=\frac{k_A}{e}\exp\left[\left(\log\left(\psi_{\lambda,k_A}\left(\frac{\tau}{\beta}\right)\right)+1\right)^{\delta_1}\right]\,,
\end{equation}
where $\delta_1=((\sigma+\tau)/\tau)^{-\lambda k_A}$. From (\ref{eq_sessanta}), (\ref{eq_espressione_beta}) and (\ref{eq_sessanta4}) one obtains
\begin{multline}\label{eq_con_log_esp_prima}
\sup_{z\in[0,\bar{\sigma}/2]}\Vert u(z,\cdot)\Vert^2_{L^2}\le C k_A \psi_{\lambda,k_A}\left(\frac{1}{\Lambda\left(\frac{1}{\tau}\log\Vert u(0,\cdot)\Vert_{L^2}\right)}\right)\times\\[2mm]
\times \exp\left\{-\frac{\sigma k_A}{e}\exp\left[\left(\log\left(\psi_{\lambda,k_A}\left(\frac{1}{\Lambda^{-1}\left(\frac{1}{\tau}\log \Vert u(0,\cdot)\Vert_{L^2}\right)}\right)\right)+1\right)^{\delta_1}\right]\right\}\times\\[2mm]
\times\left(1+\Vert u(\sigma,\cdot)\Vert^2_{L^2}\right)\,.
\end{multline}
Consider, now, the function $F$ defined by
$$
F(\zeta)\triangleq (1+D)C k_A \zeta \exp\left\{-\frac{\sigma k_A}{2e}\exp\left[\left(\log\zeta+1\right)^{\delta_1}\right]\right\}
$$
and note that
$$
\lim_{\zeta\to+\infty}F(\zeta)=0\,.
$$
Indeed, let $\epsilon>0$. It is easy to check that
\begin{multline*}
F(\zeta)<\epsilon \Leftrightarrow\exp\left\{-\frac{\sigma k_A}{2e}\exp\left[\left(\log \zeta+1\right)^{\delta_1}\right]\right\}\le \frac{\epsilon \zeta^{-1}}{C k_A}\Leftrightarrow\\[2mm]
\Leftrightarrow -\frac{\sigma k_A}{2e}\exp\left[\left(\log \zeta+1\right)^{\delta_1}\right]\le -\log\zeta+\log\frac{\epsilon}{C k_A(1+D)}\Leftrightarrow\\[2mm]
\Leftrightarrow\frac{\sigma k_A}{2e}\exp\left[\left(\log\zeta+1\right)^{\delta_1}\right]\ge \log \zeta-\log\frac{\epsilon}{Ck_A(1+D)}\Leftrightarrow\\[2mm]
\Leftrightarrow \exp\left[\left(\log\zeta+1\right)^{\delta_1}\right]\ge \frac{2e}{\sigma k_A}\log\zeta-\frac{2e}{\sigma k_A}\log\frac{\epsilon}{C k_A(1+D)}\Leftrightarrow\\[2mm]
\Leftrightarrow\left(\log\zeta+1\right)^{\delta_1}\ge \log\left(\frac{2e}{\sigma k_A}\log\zeta-\frac{2e}{\sigma k_A}\log\frac{\epsilon}{C k_A(1+D)}\right)\,,
\end{multline*}
which is true for sufficiently large $\zeta$. Analogousy, for sufficiently small $\Vert u(0,\cdot)\Vert_{L^2}$, one has
\begin{multline*}
(1+D)C k_A\psi_{\lambda,k_A}\left(\frac{1}{\Lambda^{-1}\left(\frac{1}{\tau}\log\Vert u(0,\cdot)\Vert_{L^2}\right)}\right)\times\\[2mm]
\times\exp\left\{-\frac{\sigma k_A}{2e}\exp\left[\left(\log\left(\psi_{\lambda,k_A}\left(\frac{1}{\Lambda^{-1}\left(\frac{1}{\tau}\log \Vert u(0,\cdot)\Vert_{L^2}\right)}\right)\right)+1\right)^{\delta_1}\right]\right\}\le 1\,.
\end{multline*}
So, if $\Vert u(0,\cdot)\Vert_{L^2}\le \widetilde{\rho}$, one has
\begin{equation}\label{eq_sessantasei}
\sup_{z\in[0,\bar{\sigma}/2]}\Vert u(z,\cdot)\Vert^2_{L^2}\le \exp\left\{-\frac{\sigma k_A}{2e}\exp\left[\left(\log\left(\psi_{\lambda,k_A}\left(\frac{1}{\Lambda^{-1}\left(\frac{1}{\tau}\log\Vert u(0,\cdot)\Vert\right)}\right)\right)+1\right)^{\delta_1}\right]\right\}\,.
\end{equation}
Now, since
$$
\lim_{\zeta\to +\infty}\frac{\psi_{\lambda,k_A}\left(\frac{1}{\tau}\right)}{\vert\Lambda(\zeta)\vert}=+\infty
$$
(see Lemma~\ref{lemma_per_il_logaritmo} in Appendix) for $\Vert u(0,\cdot)\Vert_{L^2}$ sufficiently small one has
$$
\psi_{\lambda,k_A}\left(\frac{1}{\Lambda^{-1}\left(\frac{1}{\tau}\left(\log\Vert u(0,\cdot)\Vert_{L^2}\right)\right)}\right)\ge\frac{1}{\tau}\left\vert \log\Vert u(0,\cdot)\Vert_{L^2}\right\vert\,.
$$
As a consequence, (\ref{eq_sessantasei}) yields
\begin{equation}\label{eq_sessantasette}
\sup_{z\in[0,\bar{\sigma}/2]}\Vert u(z,\cdot)\Vert_{L^2}^2\le \exp\left\{-\frac{\sigma k_A}{2e}\exp\left[\left(\log\left(\frac{1}{\tau}\left\vert \log\Vert u(0,\cdot)\Vert_{L^2}\right\vert\right)\right)^{\delta_1}\right]\right\}
\end{equation}
which may also be rewritten as
\begin{equation}\label{eq_sessantotto}
\sup_{z\in[0,\bar{\sigma}_1]}\Vert u(z,\cdot)\Vert_{L^2}^2\le \exp\left\{-\frac{\sigma_1 k_A}{2e}\exp\left[\left(\log\left(\frac{1}{2\tau_1}\left\vert \log\Vert u(0,\cdot)\Vert_{L^2}^2\right\vert\right)\right)^{\delta_1}\right]\right\}\,,
\end{equation}
where $\bar{\sigma}_1=\sigma_1/16$. Now, choose
$$
\tau_1=\min\left\{\frac{\sigma_1}{4},\frac{\sigma_1 k_A}{4e}\right\}
$$
and iterate the above arguments on $[\bar{\sigma}_1,T]$, finding
\begin{equation}\label{eq_sessantanove}
\sup_{z\in[\bar{\sigma}_1,\bar{\sigma}_2]}\Vert u(z,\cdot)\Vert_{L^2}^2\le \exp\left\{-\frac{\sigma_2 k_A}{2e}\exp\left[\left(\log\left(\frac{1}{2\tau_2}\left\vert \log\Vert u(\bar{\sigma}_1,\cdot)\Vert_{L^2}^2\right\vert\right)\right)^{\delta_2}\right]\right\}\,,
\end{equation}
where $\bar{\sigma}_2=\sigma_2/16$ and $\tau_2=\min\left\{\frac{\sigma}{4},\frac{\sigma_2 k_A}{4e}\right\}$. Note that
$$
\sigma_1=\min\left\{\frac{1}{\alpha_1},T^{\prime\prime}\right\}\,,\quad \sigma_2=\min\left\{\frac{1}{\alpha_1},T^{\prime\prime}-\sigma_1\right\}\,;
$$
hence $\sigma_2\le\sigma_1$ and $\tau_2\le\tau_1$. As a consequence,
\small
\begin{multline*}
\sup_{z\in[\bar{\sigma}_1,\bar{\sigma}_2]}\Vert u(z,\cdot)\Vert_{L^2}^2\le\\[2mm]
\exp\left\{-\frac{\sigma_2 k_A}{2e}\exp\left[\left(\log\left\vert\frac{1}{2\tau_2}\log\left(\exp\left\{-\frac{\sigma_1 k_A}{2e}\exp\left[\left(\log\frac{\left\vert\log\Vert u(0,\cdot)\Vert^2_{L^2}\right\vert}{2\tau_1}\right)^{\delta_1}\right]\right\}\right)\right\vert\right)^{\delta_2}\right]\right\}=\\[2mm]
=\exp\left\{-\frac{\sigma_2 k_A}{2e}\exp\left[\left(\log\left\vert-\frac{\sigma_1 k_A}{4e\tau_2}\exp\left[\left(\log\frac{\left\vert\log\Vert u(0,\cdot)\Vert^2_{L^2}\right\vert}{2\tau_1}\right)^{\delta_1}\right]\right\vert\right)^{\delta_2}\right]\right\}=\\[2mm]
=\exp\left\{-\frac{\sigma_2 k_A}{2e}\exp\left[\left(\log\frac{\sigma_1 k_A}{4e\tau_2}+\left(\log\frac{\left\vert\log\Vert u(0,\cdot)\Vert^2_{L^2}\right\vert}{2\tau_1}\right)^{\delta_1}\right)^{\delta_2}\right]\right\}\le\\[2mm]
\le \exp\left\{-\frac{\sigma_2 k_A}{2e}\exp\left[\left(\log\frac{1}{2\tau_1}\left\vert \log\Vert u(0,\cdot)\Vert^2_{L^2}\right\vert\right)^{\delta_1\delta_2}\right]\right\}\,,
\end{multline*}
\normalsize
where the last inequality holds since $\sigma_1 k_A\ge 4e\tau_2$. Merging the estimates obtained for the two intervals, yields
$$
\sup_{[0,\bar{\sigma}_2]}\Vert u(z,\cdot)\Vert^2_{L^2}\le \exp\left\{-\frac{\sigma_2 k_A}{2e}\exp\left[\log\frac{1}{2\tau_1}\left\vert\log\Vert u(0,\cdot)\Vert^2_{L^2}\right\vert\right]\right\}\,,
$$
which has the same form of the inequality obtained in $[0,\bar{\sigma}_1]$. Hence, if $T^{\prime\prime}$ is such that $0<T^\prime<T^{\prime}{\prime}<T$, iterating a finite number of times one obtains an estimate on $[0,T^\prime]$ of the form
$$
\sup_{[0,T^\prime]}\Vert u(z,\cdot)\Vert^2_{L^2}\le\exp\left\{-\widetilde{\sigma}\exp\left[\left(\log\frac{1}{2\tilde{\tau}}\left\vert\log\Vert u(0,\cdot)\Vert^2_{L^2}\right\vert\right)^{\widetilde{\delta}}\right]\right\}\,.
$$

\section*{Appendix}
\begin{lemma}\label{lemma_per_il_logaritmo}
The functions $\psi_{\lambda, k_A}$ (equation (\ref{eq_def_psi})) and $\Lambda$ (equation (\ref{eq_def_Lambda})) are such that
$$
\lim_{\zeta\to+\infty}\frac{\psi_{\lambda, k_A}\left(\frac{1}{\zeta}\right)}{\vert\Lambda(\zeta)\vert}=+\infty\,.
$$
{\bf Proof. }Note that
\begin{multline*}
\lim_{\zeta\to+\infty}\frac{\psi_{\lambda, k_A}\left(\frac{1}{\zeta}\right)}{\vert\Lambda(\zeta)\vert} = \lim_{\rho\to 0}-\frac{\rho\psi_{\lambda,k_A}(\rho)}{\phi_{\lambda,k_A}(\rho)}=\lim_{\rho\to 0}-\frac{\psi_{\lambda,k_A}(\rho)+\rho\psi^\prime_{\lambda,k_A(\rho)}}{k_A \psi_{\lambda,k_A}(\rho)}=\\[2mm]
=-\frac{1}{k_A}-\lim_{\rho\to 0}\frac{\rho \psi_{\lambda, k_A}^\prime(\rho)}{k_A \psi_{\lambda, k_A}(\rho)}=-\frac{1}{k_A}+\lim_{\rho\to 0}\frac{1}{k_A}\lambda \psi_{\lambda, k_A}(\rho)\omega\left(\frac{1}{\psi_{\lambda, k_A}(\rho)}\right)=\\[2mm]
=-\frac{1}{k_A}+\lim_{q\to 0}\frac{\lambda}{k_A}\frac{\omega(q)}{q}=-1+\lim_{q\to 0}(1-\log q)\log(1-\log q)=+\infty\,.
\end{multline*}
\end{lemma}

\bibliographystyle{plain}
\bibliography{Osgood_biblio}

\begin{thebibliography}{10}

\bibitem{mia_tesi}
D.~Casagrande.
\newblock {\em Stabilit\`a condizionata per equazioni paraboliche retrograde
  con coefficienti continui non lipsciziani}.
\newblock Master thesis (in Italian), Universit\`a degli Studi di Trieste,
  Italy, 2017, available at:
  \texttt{http:$\backslash\backslash$www.diegm.uniud.it/casagrande/tesi\_matem%
atica}.

\bibitem{NONLINAL}
D.~{Del Santo}, C.~J\"ah, and M.~Prizzi.
\newblock Conditional stability for backward parabolic equations with
  {L}og{L}ip$_t\times$ {L}ip$_x$-coefficients.
\newblock {\em Nonlinear Analysis}, 121:101--122, 2015.

\bibitem{JMPA}
D.~{Del Santo} and M.~Prizzi.
\newblock Backward uniqueness for parabolic operators whose coefficients are
  non-{L}ipschitz continuous in time.
\newblock {\em Journal de Math\'ematiques Pures et Appliqu\'ees}, 84:471--491,
  2005.

\bibitem{MATAN}
D.~{Del Santo} and M.~Prizzi.
\newblock Continuous dependence for backward parabolic operators with
  {L}og-{L}ipschitz coefficients.
\newblock {\em Mathematische Annalen}, 345:213--243, 2009.

\bibitem{AMPA}
D.~{Del Santo} and M.~Prizzi.
\newblock A new result on backward uniqueness for parabolic operators.
\newblock {\em Annali di Matematica Pura e Applicata}, 194:387--403, 2015.

\bibitem{Had_53}
J.~Hadamard.
\newblock {\em Lectures on {C}auchy's {P}roblem in {L}inear {P}artial
  {D}ifferential {E}quations}.
\newblock Yale University Press, New Haven, 1923.

\bibitem{Had_64}
J.~Hadamard.
\newblock {\em La Th\'eorie des {\'E}quations aux {D}\'eriv\'es {P}artielles}.
\newblock \'Editions Scientifique, Peking. Gauthier-Villars {\'E}diteur, Paris,
  1964.

\bibitem{Hua_Rod}
C.~Hua and L.~Rodino.
\newblock Paradifferential calculus in gevrey classes.
\newblock {\em Journal of Mathematics of Kyoto University}, 41(1):1--31, 2001.

\bibitem{Hurd}
A.E. Hurd.
\newblock Backward continuous dependence for mixed parabolic problems.
\newblock {\em Duke Mathematical Journal}, 34:493--500, 1967.

\bibitem{John}
F.~John.
\newblock Continuous dependence on data for solutions of partial differential
  equations with a prescribed bound.
\newblock {\em Comm. Pure Appl. Math.}, 13:551--585, 1960.

\bibitem{Lio_Mal}
J.-L. Lions and B.~Malgrange.
\newblock Sur l'unicit\'e r\'etrograde dans les probl\`emes mixtes
  paraboliques.
\newblock {\em Mathematica Scandinavica}, 8:277--286, 1960.

\bibitem{Pli}
A.~Pli\'s.
\newblock On non-uniqueness in {C}auchy problem for an elliptic second order
  differential equation.
\newblock {\em Bulletin of the Polish Academy of Science}, 11:95--100, 1963.

\end{thebibliography}

\end{document}